\documentclass[12pt]{amsart}
\usepackage{amssymb}
\usepackage{verbatim}
\usepackage[usenames]{color}
\usepackage{hyperref}
\usepackage{url}
\usepackage{tikz}

\newtheorem{thm}{Theorem}

\newtheorem{la}[thm]{Lemma}
\newtheorem{cor}[thm]{Corollary}

\theoremstyle{definition}
\newtheorem{df}[thm]{Definition}
\newtheorem{conv}[thm]{Convention}
\newtheorem{notat}[thm]{Notation}
\newtheorem{que}[thm]{Question}

\theoremstyle{remark}

\newtheorem{rmk}[thm]{Remark}

\newenvironment{ls}{\begin{itemize}}{\end{itemize}}
\newenvironment{lsnum}{\begin{enumerate}}{\end{enumerate}}
\newenvironment{pf}{\begin{proof}}{\end{proof}}
\newcommand{\ger}[1]{\ensuremath{\mathfrak {#1}}}
\newcommand{\scr}[1]{\ensuremath{\mathcal {#1}}}
\newcommand{\bld}[1]{\ensuremath{\mathbf {#1}}}
\newcommand{\bbb}[1]{\ensuremath{\mathbb {#1}}}

\newcommand{\dom}[1]{\ensuremath{{\text{Dom}}(#1)}}
\newcommand{\ran}[1]{\ensuremath{{\text{Range}}(#1)}}

\newcommand{\emp}{\varnothing}
\renewcommand{\phi}{\varphi}
\newcommand{\eps}{\varepsilon}
\newcommand{\sq}[1]{\ensuremath{\langle#1\rangle}}

\newcommand{\restr}{\mathop{\upharpoonright}}

\newcommand{\notarrow}{\kern .42em\not\kern -.42em\longrightarrow}
\renewcommand{\th}{\ensuremath{{}^{\text{th}}}}

\newcommand{\rr}{\ensuremath{\ger{rr}}}
\newcommand{\rrf}{\ensuremath{\ger{rr}_f}}
\newcommand{\rri}{\ensuremath{\ger{rr}_i}}
\newcommand{\rro}{\ensuremath{\ger{rr}_o}}
\newcommand{\cov}[1]{\ensuremath{\bld{cov}(\scr{#1})}}
\newcommand{\non}[1]{\ensuremath{\bld{non}(\scr{#1})}}
\newcommand{\nm}[1]{\ensuremath{\left\Vert #1 \right\Vert}}
\newcommand{\nmb}[1]{\ensuremath{\big\Vert #1 \big\Vert}}
\newcommand{\masc}{\ensuremath{\text{MA}(\sigma\text{-centered})}}

\newcommand{\noprint}[1]{\relax}

\title{The Rearrangement Number}
\author[Blass]{Andreas Blass}
\address[A.~R.~Blass]{Mathematics Department\\
University of Michigan\\
Ann Arbor, MI 48109--1043, U.S.A.}
\email{ablass@umich.edu}
\urladdr{http://www.math.lsa.umich.edu/~ablass/}
\author[Brendle]{J\"org Brendle}
\address[J.~Brendle]{Graduate School of System Informatics, Kobe
  University, 1--1   Rokkodai, Nada-ku, 657-8501 Kobe,  Japan}
\email{brendle@kobe-u.ac.jp}
\author[Brian]{Will Brian}
\address[W.~Brian]{Department of Mathematics, Baylor University, One Bear Place
  \#97328, Waco, TX 76798-7328, U.S.A.}
\email{wbrian.math@gmail.com}
\urladdr{http://wrbrian.wordpress.com}
\author[Hamkins]{Joel David Hamkins}
\address[J.~D.~Hamkins]{Mathematics, The Graduate Center of the City
  Univeristy of New York, 365 Fifth Avenue, New York, NY 10016,
  U.S.A. and Mathematics, College of Staten Island of CUNY, Staten
  Island, NY 10314, U.S.A.}
\email{jhamkins@gc.cuny.edu}
\urladdr{http://jdh.hamkins.org}
\author[Hardy]{Michael Hardy}
\address[M.~Hardy]{Department of Mathematics, Hamline University,
  Saint Paul, MN 55104, U.S.A.}
\email{drmichaelhardy@gmail.com}
\author[Larson]{Paul B. Larson}
\address[P.~B.~Larson]{Department of Mathematics, Miami University,
  Oxford, OH 45056, U.S.A.}
\email{larsonpb@miamioh.edu}
\urladdr{http://www.users.miamioh.edu/larsonpb/}

\thanks{The research of the first and fourth authors on this topic
  took place in part while they were research fellows at the Isaac Newton Institute
  for Mathematical Sciences in the program ``Mathematical,
  Foundational and Computational Aspects of the Higher Infinite
  (HIF).''  Most of the first author's writing was done while he was a
  visiting scientist at the Simons Institute for Theory of Computing
  in Berkeley. The research of the second author was partially
  supported by Grant-in-Aid for Scientific Research (C) 15K04977,
  Japan Society for the Promotion of Science.  The research of the
  fourth author was supported in part by Simons Foundation grant
  209252.  The research of the sixth author was supported in part by
  NSF grant DMS-1201494. The fourth author thanks Heike Mildenberger
  for discussions about our topic.  We are also grateful for
  additional answers to \cite{HardyQ} posted on MathOverflow by Robert
  Israel and Aaron Meyerowitz.
   \vspace{.5mm} \\ 
  Commentary concerning this article can be made on the fourth author's
blog at 
\url{http://jdh.hamkins.org/the-rearrangment-number}. \vspace{.5mm} \\
This paper is written to be accessible to students at the graduate level and those who are not experts in set theory. A streamlined version of this paper, more suitable for experts, can be accessed at \\ 
\url{https://arxiv.org/abs/1612.07830/v2}.
}

\begin{document}

\begin{abstract}
  How many permutations of the natural numbers are needed so that
  every conditionally convergent series of real numbers can be
  rearranged to no longer converge to the same sum? We show that the
  minimum number of permutations needed for this purpose, which we
  call the rearrangement number, is uncountable, but whether it equals
  the cardinal of the continuum is independent of the usual axioms of
  set theory.  We compare the rearrangement number with several
  natural variants, for example one obtained by requiring the
  rearranged series to still converge but to a new, finite limit.  We
  also compare the rearrangement number with several well-studied
  cardinal characteristics of the continuum. We present some new
  forcing constructions designed to add permutations that rearrange
  series from the ground model in particular ways, thereby obtaining
  consistency results going beyond those that  follow from comparisons
  with familiar cardinal characteristics.  Finally we deal briefly
  with some variants concerning rearrangements by a special sort of
  permutations and with rearranging some divergent series to become
  (conditionally) convergent.
\end{abstract}

\maketitle

\section{Introduction}          \label{intro}

Let $\sum_na_n$ be a conditionally convergent series of real numbers.
Permutations $p$ of the summands, producing rearrangements
$\sum_na_{p(n)}$ of the original series, can disrupt the convergence
in several ways. Riemann \cite{riem} showed that one can obtain
rearrangements converging to any prescribed real number. (At the end
of this introduction, we record some information about the history of
Riemann's Rearrangement Theorem.)  Minor modifications of Riemann's
argument produce rearrangements that diverge to $+\infty$ or to
$-\infty$, as well as rearrangements that diverge by oscillation.

Instead of arbitrary permutations of $\bbb N$, can some limited class
$C$ of permutations suffice to disrupt, in one way or another, the
convergence of all conditionally convergent series?  In particular,
how small can the cardinality $|C|$ of such a class $C$ be? This
question was raised on MathOverflow by the fifth author of this paper
\cite{HardyQ}, and there were partial answers and comments from
several of the other authors \cite{BrianA, HamkinsA, LarsonA}.

It turns out that the question has substantial set-theoretic
ramifications, and the present paper reports on these, as well as some
related matters.  We begin, in
Section~\ref{def}, by defining the \emph{rearrangement numbers,} the
minimal cardinalities of families of permutations needed to achieve
various sorts of disruption of the convergence of conditionally
convergent series.  We point out easy connections between some of
these cardinal numbers.  This section also contains a few notational
conventions that we use throughout the paper.

In Section~\ref{pad1}, we introduce a technique of padding series by
inserting many zero terms, and we use this method to demonstrate that
the rearrangement numbers are uncountable. That is, countably many
permutations cannot suffice to disrupt all conditionally convergent
series.  The same method is applied in a more sophisticated way in
Section~\ref{pad2} to establish stronger lower bounds for the
rearrangement numbers.

In the two intervening sections, we introduce and exploit another
technique, to ``mix'' several given permutations into a single one that
accomplishes some of the disruption of the given ones. In
Section~\ref{mix1}, we use this technique to show that several of the
rearrangement numbers coincide.  In Section~\ref{mix2}, we extend the
technique to relate rearrangement to properties of Baire
category.

In Section~\ref{signs}, we relate rearrangement to properties of
Lebesgue measure, using information about series in which the signs
of the terms are chosen at random.

Sections~\ref{force1} and \ref{force2} are devoted to showing that
even the largest of the rearrangement numbers can consistently be
strictly smaller than the cardinality of the continuum.  (For the
smallest of the rearrangement numbers, this consistency already
follows from the result in Section~\ref{mix2}.)

The last two sections treat related questions in modified contexts.
Section~\ref{ad} looks at divergent series that can be rearranged to
(conditionally) converge. We show that, to achieve this sort of
rearrangement for all such series, one must use as many permutations
as the cardinality of the continuum.  Section~\ref{shuffles} is about
the situation where, as in the usual proof of Riemann's Rearrangement
Theorem, the permutations leave the relative order of the positive
terms and the relative order of the negative terms unchanged, i.e.,
they affect only the interleaving between the positive and negative
terms.

Except for the last two sections, we have tried to arrange the
material in order of increasing set-theoretic prerequisites.
Sections~\ref{def}, \ref{pad1}, and \ref{mix1}  use only
elementary cardinal arithmetic.  Sections~\ref{mix2}, \ref{pad2}, and
\ref{signs}  involve cardinal
characteristics of the continuum.  We give the relevant definitions,
but we refer to \cite{hdbk} for some facts about these
characteristics.  The results in Sections~\ref{force1} and
\ref{force2} presuppose general knowledge of forcing, including
finite-support iterations.  In the last two sections, \ref{ad} and
\ref{shuffles}, where we work in modified contexts, we return to more
elementary methods.  Most of the section headings indicate not only
the results obtained in the section but also the methods used to
obtain them.

\subsection*{History of Riemann's theorem}
On page 97 of \cite{riem}, Riemann gives the nowadays familiar proof
that a conditionally convergent series can be rearranged to converge
to any prescribed finite sum. Riemann died in 1866, and the
publication of this paper in 1867 was arranged by Dedekind.  A
footnote by Dedekind on the first page says that the paper was
submitted by Riemann for his habilitation at G\"ottingen in 1854.
Dedekind further explains that Riemann apparently didn't intend to
publish it, but that its publication is justified by its intrinsic
interest and its method of treating the basic principles of
infinitesimal analysis.  The footnote is dated 1867.  (The next paper
in the journal, also dating from 1854 and published posthumously
through the efforts of Dedekind, is the one where Riemann introduces
what is now called Riemannian geometry.)

Riemann attributes to Dirichlet the observation that there is a
crucial difference between (what are nowadays called) absolutely and
conditionally convergent series. Riemann immediately continues with
the construction of a rearrangement achieving any prescribed sum for a
conditionally convergent series.  The possibility of such
rearrangements does not seem to be in the cited paper of Dirichlet
\cite{dirichlet}.  What Dirichlet does point out there is that a
conditionally convergent series can become divergent when the terms
are (not rearranged but) multiplied by factors that approach 1.  This
is on page~158, and the mention of conditional convergence is only in a
parenthetical comment.  Specifically, Dirichlet refers to an argument
in which Cauchy inferred the convergence of a series from the
convergence of another series whose corresponding terms differ by
factors approaching 1; Dirichlet gives an example showing that this
inference is not valid ``(du moins lorsque, comme il arrive ici, les
termes n'ont pas tous le m\^eme signe).''  Incidentally, this paper
introduced what is now called the Dirichlet kernel in the theory of
Fourier series.

Finally, we briefly mention some more recent history for readers 
interested in constructive mathematics.  Diener and Lubarsky \cite{dl} have
studied a version Riemann's Rearrangement Theorem in the context of
constructive mathematics.  Since classically equivalent formulations
need not be constructively equivalent, it it appropriate to specify
that they use the formulation ``If every permutation of a series of
real numbers converges, then the series converges absolutely.''  (In
constructive logic, it seems that this formulation neither implies nor
is implied by ``If a series converges conditionally (i.e., not
absolutely), then some rearrangement diverges.'')  Diener and Lubarsky
show that this form of Riemann's Rearrangement Theorem is not
constructively provable; they exhibit a topological model where it
fails.

\section{Definitions, Conventions, and Basic Facts} \label{def}

We define the \emph{rearrangement numbers}, denoted by \rr, sometimes
with subscripts, as the minimum number of permutations of \bbb N
needed to disrupt, in various ways, the convergence of all
conditionally convergent series.

\begin{df}      \label{df-rr}
  \rr\ is the smallest cardinality of any family $C$ of
      permutations of \bbb N such that, for every conditionally
      convergent series $\sum_na_n$ of real numbers, there is some
      permutation $p\in C$ for which the rearrangement
      $\sum_na_{p(n)}$ no longer converges to the same limit.
\end{df}

In this definition, the rearranged series might converge to a (finite)
sum different from $\sum_na_n$, might diverge to $+\infty$ or to
$-\infty$, or might diverge by oscillation.  If we specify one of
these options, we get more specific rearrangement numbers, as follows. 

\vspace{15mm}

\begin{df}      \label{df-rrsub}
\hskip1pt
\begin{ls}
  \item \rrf\ is defined like \rr\ except that $\sum_na_{p(n)}$ is
    required to converge to a finite sum (different from $\sum_na_n$).
\item \rri\ is defined like \rr\ except that $\sum_na_{p(n)}$ is
    required to diverge to $+\infty$ or to $-\infty$.
\item \rro\ is defined like \rr\ except that $\sum_na_{p(n)}$ is
    required to diverge by oscillation.
\end{ls}
\end{df}

The subscripts $f$, $i$, and $o$ are intended as mnemonics for
``finite,'' ``infinite,'' and ``oscillate.''

We shall not discuss variants in which the disruption of convergence
is even more specific, for example by distinguishing oscillation
between finite bounds from oscillation over the whole real line.  We
shall, however, have use for variants in which two of the three sorts
of disruption are allowed; we denote these by \rr\ with two
subscripts.  Thus, for example, $\rr_{fi}$ is the minimum size of a
set $C$ of permutations of \bbb N such that, for every conditionally
convergent $\sum_na_n$, there is $p\in C$ for which $\sum_na_{p(n)}$
either converges to a different finite sum or diverges to $+\infty$ or
to $-\infty$.  The definitions of $\rr_{fo}$ and $\rr_{i o}$ are
analogous.  Of course, if we allow all three sorts of disruption, we
could write $\rr_{fio}$, but we have chosen (in
Definition~\ref{df-rr}) to denote this cardinal simply by \rr, because
it seems to be the most natural of all these rearrangement numbers and
because it was the subject of the original question in \cite{HardyQ}.

We have, a priori, seven rearrangement numbers: The original \rr,
three variants with one subscript, and three with two subscripts.
Clearly, if one variant allows more modes of disruption than another,
then every family of permutations adequate for the latter variant is
also adequate for the former, and therefore the former cardinal is
less than or equal to the latter.  Figure~\ref{fig1} is a Hasse
diagram showing the seven variants and the order relationships
resulting from this elementary observation.  Later, we shall see that
some of these variants coincide, so the diagram will become simpler.

\begin{center}
\begin{figure}
\begin{tikzpicture}[, xscale=2,yscale=2]

\draw (0,0) -- (-1,1);
\draw (0,0) -- (0,1);
\draw (0,0) -- (1,1);
\draw (-1,1) -- (-1,2);
\draw (-1,1) -- (0,2);
\draw (0,1) -- (-1,2);
\draw (0,1) -- (1,2);
\draw (1,1) -- (0,2);
\draw (1,1) -- (1,2);

\draw [fill=white,white] (0,0) ellipse (5pt and 3pt);  \node at (0,0) {$\rr$};
\draw [fill=white,white] (-1,1) ellipse (5pt and 4pt);  \node at (-1,1) {$\rr_{fi}$};
\draw [fill=white,white] (-.01,1.01) ellipse (5.5pt and 5pt);  \node at (0,1) {$\rr_{fo}$};
\draw [fill=white,white] (1,1) ellipse (5pt and 4pt);  \node at (1,1) {$\rr_{io}$};
\draw [fill=white,white] (-1,2) ellipse (6pt and 4pt);  \node at (-1,2) {$\rrf$};
\draw [fill=white,white] (0,2) ellipse (6pt and 4pt);  \node at (0,2) {$\rri$};
\draw [fill=white,white] (1,2) ellipse (5pt and 3pt);  \node at (1,2) {$\rro$};

\end{tikzpicture}
\caption{The seven \textit{a priori} rearrangement numbers}
\label{fig1}
\end{figure}
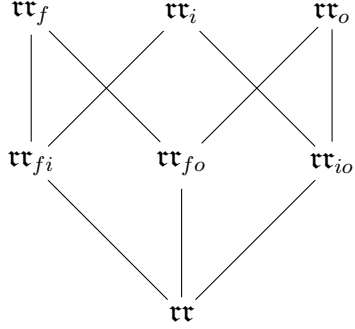
\end{center}

Before proceeding, it will be convenient to adopt the following
conventions, which are standard in set theory.

\begin{conv}    \label{conven}
Each natural number $n$ is identified with the set of strictly smaller
natural numbers, $n=\{0,1,\dots,n-1\}$.  If $f$ is a function and $S$
is a subset of its domain, then $f[S]=\ran{f\restr S}=\{f(x):x\in
S\}$.  Note that we use round parentheses for the value of a function
at a point in its domain and square brackets for the image of a subset
of the domain. In particular, if $f$ has domain \bbb N, then, for
example, $f[3]=\{f(0),f(1),f(2)\}$ and this is usually different from
$f(3)$.
\end{conv}

\begin{notat}
  The cardinality $2^{\aleph_0}$ of the continuum is denoted by \ger
  c.
\end{notat}

Riemann's Rearrangement Theorem and the minor modifications that
achieve divergence to $\pm\infty$ and divergence by oscillation
immediately imply that \ger c is an upper bound for all seven of our
rearrangement numbers.

\section{Padding with Zeros, Uncountability} \label{pad1}

In this section, we obtain our first lower bound for rearrangement
numbers.  Stronger lower bounds will be obtained later.  Recall that
$\rr$ is the smallest of the rearrangement numbers, so the following
result, stated for \rr, implies the same for all the other
rearrangement numbers.

\begin{thm}     \label{uncbl:thm}
  \rr\ is uncountable.
\end{thm}

\begin{pf}
  We must show that, given any countable set $C=\{p_n:n\in\bbb N\}$ of
  permutations of \bbb N, there is a conditionally convergent series
  $\sum_na_n$ such that, for each permutation $p\in C$, the rearranged
  series $\sum_na_{p(n)}$ converges to the same sum as the original
  $\sum_na_n$.  For this purpose, we start with any conditionally
  convergent series $\sum_nb_n$, for example the alternating harmonic
  series $\sum_n(-1)^n/n$, and we modify it by inserting a large
  number of zeros between consecutive terms.  The purpose of the zeros
  is to put the non-zero terms so far apart that the permutations in
  $C$ leave their ordering essentially unchanged.

  To make this strategy precise, we begin by defining a rapidly
  increasing function $l:\bbb N\to\bbb N$ by the following induction;
  the intention is that $l(k)$ tells the location where $b_k$ should
  go in the padded series. Begin by setting $l(0)=0$.  After $l(k)$
  has been defined, choose $l(k+1)$  larger than $l(k)$ and
  different from the finitely many numbers of the form $p_m(j)$
  for $m\leq k$ and $j\leq p_m^{-1}(l(k))$.  This definition ensures that
\[
(\forall m\leq k)\ p_m^{-1}(l(k))<p_m^{-1}(l(k+1)).
\]
That is, the relative order of $l(k)$ and $l(k+1)$ is not changed if
we apply to them the inverses of any of the first $k+1$ elements
$p_0,\dots,p_k$ of $C$.  Equivalently, the inverse of any particular
$p_m\in C$ preserves the relative order of all but the first $m$
elements of the range of $l$. (Our preoccupation here with the
inverses of the permutations from $C$, rather than with the
permutations themselves, is explained by the fact that the summand
that appears in position $k$ in a series $\sum_na_n$ gets moved to
position $p^{-1}(k)$ in the rearranged series $\sum_na_{p(n)}$.)

As indicated above, we let $l(k)$ tell us the location where the
$k\th$ term $b_k$ of our original series $\sum_nb_n$ should go in our
padded series $\sum_na_n$.  That is, we define
\[
a_n=
\begin{cases}
  b_k&\text{if }n=l(k),\\
0&\text{if } n\notin\ran l.
\end{cases}
\]
Thus, the series $\sum_na_n$ has the same nonzero terms in the same
order as $\sum_nb_n$; the only difference is that many zeros have been
inserted.  So $\sum_na_n$ is conditionally convergent (to the same sum
as $\sum_nb_n$).  Furthermore, if we apply any permutation $p_m\in C$
to $\sum_na_n$, all but the first $m$ nonzero terms will remain in
their original order.  Thus, the rearranged series $\sum_na_{p_m(n)}$
still converges to the same sum.

That is, the countable set $C$ is not as required in the definition of \rr.
\end{pf}

Theorem~\ref{uncbl:thm} and the observations in the preceding section
tell us that all our rearrangement numbers lie between $\aleph_1$ and
\ger c, inclusive.  So they qualify as cardinal characteristics of the
continuum \cite{hdbk}.  Like all such characteristics, they are
uninteresting if the continuum hypothesis (CH) holds, i.e., if
$\aleph_1=\ger c$.  But in the absence of CH, it is reasonable to ask
how they compare with more familiar cardinal characteristics like
those described in \cite{hdbk}.  We shall obtain numerous such
comparisons in the following sections.

\section{Oscillation Is Easy, Mixing Permutations} \label{mix1}

In this section, we show that our seven rearrangement numbers are in
fact only four, because several of them provably coincide.

\begin{thm}     \label{mix-o:thm}
  $\rr=\rr_{fo}=\rr_{io}=\rro$.
\end{thm}

\begin{pf}
  In view of the observations in Section~\ref{def} (see
  Figure~\ref{fig1}), it suffices to prove $\rro\leq\rr$.  The key
  ingredient in the proof is the following lemma, which shows how to
  ``mix'' two permutations.  Recall from Convention~\ref{conven} that
  $g[n]$ means $\{g(0),g(1),\dots,g(n-1)\}$.

\begin{la}      \label{mix:la}
  For any permutation $p$ of \bbb N, there exists a permutation $g$ of
  \bbb N such that
  \begin{ls}
    \item $g[n]=p[n]$ for infinitely many $n$,
      and
\item $g[n]=n$ for infinitely many $n$
  \end{ls}
\end{la}

\begin{pf}[Proof of Lemma]
  Let an arbitrary permutation $p$ of \bbb N be given.  Notice that,
  for any $n$ and any injective function $h:n\to\bbb N$, there is some
  $M>n$ such that $h$ can be extended to an injective function
  $h':M\to\bbb N$ with $\ran{h'}=p[M]$.  Indeed, since $p$
  maps onto \bbb N, we can choose $M$ so large that
  $h[n]\subseteq p[M]$, and then we can extend
  $h$ to the required $h'$ by sending the elements of $M\setminus n$
  bijectively to the $M-n$ elements of
  $p[M]\setminus h[n]$.

Similarly, using the identity function on \bbb N instead of $p$, we
can extend any injective $h:n\to\bbb N$ to an injective $h':M\to\bbb
N$ such that $\ran{h'}=M$.

Applying these two constructions, one for $p$ and one for the identity
map, alternately, we obtain the desired $g$.  That is, we define, by
induction on $k$, injective functions $g_k:n_k\to\bbb N$ to serve as
initial segments of $g$. We start with the empty function as $g_0$. If
$g_k$ is already defined for an odd number $k$, then we properly
extend it to some $g_{k+1}:n_{k+1}\to\bbb N$ such that
$\ran{g_{k+1}}=p[n_{k+1}]$.  If $g_k$ is already defined for an even
number $k$, then we properly extend it to some
$g_{k+1}:n_{k+1}\to\bbb N$ such that $\ran{g_{k+1}}=n_{k+1}$. Since
each $g_{k+1}$ is injective and extends $g_k$, and since every natural
number eventually appears in the range of some $g_k$, there is a
permutation $g$ of \bbb N that extends all of the $g_k$'s.  Our
construction for odd $k$ ensures that $g$ satisfies the first
conclusion of the lemma, and our construction for even $k$ ensures
that it satisfies the second. This completes the proof of the lemma.
\end{pf}

Returning to the proof of the theorem, suppose $C$ is a set of
permutations of \bbb N as in the definition of \rr. For each $p\in C$,
let $g_p$ be a permutation of \bbb N that mixes $p$ and the identity,
as in the lemma.  Define $C'=C\cup\{g_p:p\in C\}$.  Because $C$ is
infinite (in fact uncountable, by Theorem~\ref{uncbl:thm}), $C$ and
$C'$ have the same cardinality, so the theorem will be proved if we
show that $C'$ is as in the definition of \rro.

Let $\sum_na_n$ be an arbitrary conditionally convergent series of
real numbers; we must find a permutation in $C'$ that makes the series
diverge by oscillation.  There is a permutation $p\in C$ that disrupts
the convergence of $\sum_na_n$. If $\sum_na_{p(n)}$ diverges by
oscillation, then we are done, since $p\in C'$.

Suppose $\sum_na_{p(n)}$ converges to a finite sum $t$ different from the
sum $s$ of $\sum_na_n$.  Then, by our choice of $g_p$, infinitely
many of the partial sums of $\sum_na_{g_p(n)}$ agree with the
corresponding partial sums of $\sum_na_{p(n)}$ and thus approach $t$,
while infinitely many other partial sums of $\sum_na_{g_p(n)}$ agree
with those of $\sum_na_n$ and thus approach $s$.  Since $s\neq t$, we
conclude that $\sum_na_{g_p(n)}$ diverges by oscillation.

Finally, if $\sum_na_{p(n)}$ diverges to $+\infty$ or to $-\infty$, the
same argument as in the preceding paragraph again shows that
$\sum_na_{g_p(n)}$ diverges by oscillation.
\end{pf}

In view of Theorem~\ref{mix-o:thm}, the Hasse diagram
of rearrangement
numbers in Figure~\ref{fig1} simplifies to Figure~\ref{fig2}.

\begin{center}
\begin{figure}
\begin{tikzpicture}[, xscale=2,yscale=2]

\draw (0,0) -- (0,1);
\draw (0,1) -- (-1,2);
\draw (0,1) -- (1,2);

\draw [fill=white,white] (0,0) ellipse (5.5pt and 3.5pt);  \node at (0,0) {$\rr$};
\draw [fill=white,white] (-.01,1.01) ellipse (6.5pt and 5.5pt);  \node at (0,1) {$\rr_{fi}$};
\draw [fill=white,white] (-1,2) ellipse (6.5pt and 4.5pt);  \node at (-1,2) {$\rrf$};
\draw [fill=white,white] (1,2) ellipse (5.5pt and 3.5pt);  \node at (1,2) {$\rri$};

\end{tikzpicture}
  \caption{The four rearrangement numbers}
\label{fig2}
\end{figure}
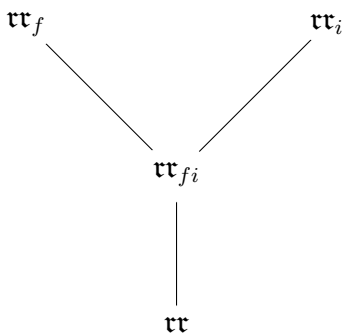
\end{center}

\section{Baire Category, More Mixing}   \label{mix2}

This is the first of several sections in which certain well-known cardinal
characteristics of the continuum play a role.  We define these
characteristics here (even though some will be needed only in later
sections), and we list the known inequalities relating them.  See
\cite{hdbk} for more information about these (and other) cardinal
characteristics.

We begin with two characteristics related to the ordering of sequences
of natural numbers by eventual domination.

\begin{df}      \label{bd:df}
For functions $f,g:\bbb N\to\bbb N$, define $f\leq^*g$ to mean that
$f(n)\leq g(n)$ for all but finitely many $n\in\bbb N$.  The
\emph{bounding number} \ger b is the minimum cardinality of a family
\scr B of functions $f:\bbb N\to\bbb N$ such that no single $g$ is
$\geq^*$ all $f\in\scr B$.  The \emph{dominating number} is the
minimum cardinality of a family \scr D of functions $f:\bbb N\to\bbb
N$ such that every $g:\bbb N\to\bbb N$ is $\leq^*$ at least one member
of \scr D.
\end{df}

Next are two characteristics arising from the Baire Category Theorem.

\begin{df}      \label{baire:df}
  A subset $M$ of a complete metric space $X$ is \emph{meager} (also
  called \emph{first category}) if it can be covered by countably many
  closed sets with empty interiors in $X$.  A \emph{comeager} set is
  the complement of a meager set; equivalently, it is a set that
  includes the intersection of countably many dense open subsets of
  $X$.  When $X$ is the space $\bbb N^{\bbb N}$ of functions $\bbb
  N\to\bbb N$, equipped with the product topology, we denote the
  family of meager subsets of $\bbb N^{\bbb N}$ by \scr M.  The
  \emph{covering number for Baire category}, \cov M, is the minimum
  number of meager sets needed to cover $\bbb N^{\bbb
    N}$. Equivalently, it is the minimum number of dense open subsets
  of $\bbb N^{\bbb N}$ with empty intersection.  The \emph{uniformity
    of Baire category}, \non M, is the minimum cardinality of a
  non-meager subset of $\bbb N^{\bbb N}$.
\end{df}

Both \cov M and \non M would be unchanged if we used the real line in
place of $\bbb N^{\bbb N}$ in the definition.

These characteristics for Baire category have analogs for Lebesgue
measure.

\begin{df}      \label{lebesgue:df}
By \emph{Lebesgue measure} we mean the product measure on the set
$2^{\bbb N}$ of functions $\bbb N\to2$ induced by the uniform measure
on $2$.  This is equivalent to the usual Lebesgue measure on the real
interval [0,1] via binary expansions of reals. We write \scr N for the
family of subsets of measure zero in $2^{\bbb N}$.  The \emph{covering
  number for Lebesgue measure}, \cov N, is the minimum number
of measure-zero sets needed to cover $2^{\bbb N}$.  The
\emph{uniformity of Lebesgue measure}, \non N, is the smallest
cardinality of a subset of $2^{\bbb N}$ that does not have measure
zero (i.e., that has positive outer measure).
\end{df}

All six of the cardinal characteristics defined here, \ger b, \ger d,
\cov M, \non M, \cov N, and \non N, lie between $\aleph_1$ and \ger c,
inclusive.  There are several known inequalities between them:
Obviously $\ger b\leq\ger d$.  An analysis of the nature of meager
sets in $\bbb N^{\bbb N}$ shows that $\ger b\leq\non M$ and
$\cov M\leq\ger d$.  Finally, a result of Rothberger \cite{roth}
(\cite[Theorem~5.11]{hdbk} is a more accessible reference) asserts
that $\cov N\leq\non M$ and $\cov M\leq\non N$.  No further
inequalities between these six cardinals are provable in ZFC.  In
fact, given any assignment of values $\aleph_1$ or $\aleph_2$ to these
cardinals and to \ger c, if the assignment is consistent with the
inequalities stated here then it is realized in some models of ZFC;
the relevant models can be found in \cite[Chapter~7]{bj}.

The main result of this section is our only nontrivial (i.e., other
than \ger c) upper bound on a rearrangement number.  Unfortunately, it
applies only to the smallest of the four rearrangement numbers, \rr.
We note that the main idea of the following proof, which is to examine the effect of ``generic'' permutations on a conditionally convergent series, has been considered before, e.g. in \cite{agnew} or in \cite[p. 346]{cohen}. 

\begin{thm}     \label{belownonM:thm}
$\rr\leq\non M$.
\end{thm}

\begin{pf}
  Although \non M was defined using the space $\bbb N^{\bbb N}$ of all
  functions $\bbb N\to\bbb N$, we shall use the fact that \non M would
  be unchanged if we used instead the subspace consisting of only the
  permutations of \bbb N.  In fact, the subspace is homeomorphic to
  the whole space.  This follows, via \cite[Theorem~7.7]{kechris},
  from the easily verified facts that the subspace is a $G_\delta$ set
  in $\bbb N^{\bbb N}$ and that it has no nonempty compact open
  subset.  A direct construction of a homeomorphism can be given by
  coding permutations as follows. Think of a permutation $p$ as a set
  of ordered pairs $(n,p(n))$, and list these ordered pairs in a
  sequence \sq{(a_k,b_k):k\in\bbb N} in such a way that, for each
  even $k$, $a_k$ is the smallest number that has not occurred as
  $a_j$ for any $j<k$, while for each odd $k$, $b_k$ is the smallest
  number that has not occurred as $b_j$ for any $j<k$. For each even
  $k$, consider where $b_k$ occurs in the increasing enumeration of
  all the natural numbers different from $b_j$ for all $j<k$; let
  $c_k$ be the number of this position.  Similarly, for each odd $k$,
  consider where $a_k$ occurs in the increasing enumeration of all the
  natural numbers different from $a_j$ for all $j<k$, and let $c_k$ be
  the number of this position.  Then the correspondence between
  permutations $p$ and the sequences
  $\bar c=\sq{c_k:k\in\bbb N}\in\bbb N^{\bbb N}$ defined in this
  manner is easily seen to be a bijection.  Since any finite part of
  $\bar c$ is determined by a finite part of $p$ and vice versa, this
  bijection is continuous in both directions, i.e., it is a
  homeomorphism.

Thus, to prove the theorem, it suffices to show that any
  nonmeager set $C$ of permutations suffices to disrupt the
  convergence of all conditionally convergent series.

Consider any conditionally convergent series $\sum_na_n$. We claim
that the permutations $p$ such that $\sum_na_{p(n)}$ fails to converge
to the same sum as $\sum_na_n$ form a comeager set.  Since $C$ isn't
meager, it must meet this comeager set, and that suffices to complete
the proof of the theorem.

In fact, we prove a stronger claim, namely that the set of
permutations $p$ such that a subsequence of the partial sums of
$\sum_na_{p(n)}$ diverges to $+\infty$ is the intersection of
countably many dense open subsets of the space of all permutations.
Indeed, this set is
\[
\bigcap_{k\in\bbb N}\bigcup_{m\geq k}\{p:\sum_{n=0}^ma_{p(n)}\geq k\},
\]
so we need only prove that, for each $k$, the set
\[
U_k=\bigcup_{m\geq k}\{p:\sum_{n=0}^ma_{p(n)}\geq k\}
\]
is dense and
open.  Taking into account the definition of the topology on the set
of permutations, as a subspace of the product space $\bbb N^{\bbb N}$,
we see that openness is immediate: If $p\in U_k$, then there is
$m\geq k$ with $\sum_{n=0}^ma_{p(n)}\geq k$, and any $p'$ that agrees
with $p$ up to $m$ is also in $U_k$.

It remains to show that $U_k$ is dense and, in view of the definition
of the topology, what we must show is that every injective function
$h:d\to\bbb N$, for any $d\in\bbb N$, can be extended to an injective
function $h':d'\to\bbb N$ with $\sum_{n=0}^ma_{h'(n)}\geq k$ for some
$m\geq k$.  But this is easy. First extend $h$, if necessary, so that
its domain is $\geq k$.  Then extend it further, using successive
values $h'(n)$ for which $a_{h'(n)}$ is positive, until the sum
exceeds $k$.  The sum will eventually exceed $k$ because, in a
conditionally convergent series like $\sum_na_n$, the sum of the
positive terms diverges to $+\infty$.
\end{pf}

A slight modification of this proof establishes the comeagerness of
the set of permutations $p$ for which $\sum_na_{p(n)}$ has not only
arbitrarily high positive partial sums but also arbitrarily low
negative partial sums.  Thus, any nonmeager set of permutations
suffices to convert any conditionally convergent series to a series
that diverges by oscillation over the whole real line, i.e., has
partial sums ranging from $-\infty$ to $+\infty$ and thus (since the
terms of the series approach zero) dense in the real line.

The argument in the proof cannot, however, be modified to obtain
convergence of $\sum_na_{p(n)}$ to a different finite sum or divergence
to $+\infty$ or $-\infty$.  That is, the argument does not make \non M
an upper bound for \rrf\ or \rri\ or even $\rr_{fi}$.  We shall see
later that this is not a defect of the argument; these larger
rearrangement numbers can consistently be larger than \non M.

To close this section, we establish a Baire category lower bound for
$\rr_{fi}$.  Even though this bound will be superseded by a stronger
one in the next section, the argument seems to be of sufficient
interest to justify mentioning it here.

\begin{thm}  \label{abovecovM:thm}
$\rr_{fi}\geq\cov M$.
\end{thm}

\begin{pf}
We begin by looking more closely at the proof of
Lemma~\ref{mix:la}. We showed there that, given a permutation $p$ of
\bbb N, we can extend any injective $h:n\to\bbb N$ to an injective
function $h': M\to\bbb N$ such that $\ran{h'}=p[M]$.  From this it
follows that the set $A_p$ defined as
\[
\{g: g\text{ a permutation of }\bbb N\text{ and }g[M]=p[M]
\text{ for infinitely many }M\}
\]
is a comeager subset of the space of all permutations of \bbb N.

Now suppose, toward a contradiction, that $\rr_{fi}<\cov M$, and let
$C$ witness this.  That is, $|C|<\cov M$ and every conditionally
convergent series can be rearranged by a permutation in $C$ so as to
either converge to a different sum or diverge to $+\infty$ or
$-\infty$.  Without loss of generality, let the identity permutation
id be a member of $C$.  The set $\bigcap_{p\in C}A_p$, being the
intersection of fewer than \cov M comeager sets, is nonempty, so let
$g$ be a member of it.

Consider any conditionally convergent series $\sum_na_n$ and, by our
choice of $C$, let $p\in C$ be such that $\sum_na_{p(n)}$ has either a
different finite sum from $\sum_na_n$  or an infinite sum.  Since
$g\in A_p$, we have $g[M]=p[M]$ for infinitely many $M$ and also
$g[M]=M$ for another infinitely many $M$.  As a result, the same
argument as in the proof of Theorem~\ref{mix-o:thm} shows that
$\sum_na_{g(n)}$ diverges by oscillation.

Thus, $\{g\}$ witnesses that $\rro=1$, which is absurd, by
Theorem~\ref{uncbl:thm}.  This contradiction completes the proof that
$\rr_{fi}\geq\cov M$.
\end{pf}

It follows from the preceding results that the rearrangment numbers
are not all provably equal.

\begin{cor}     \label{separaterr:cor}
It is consistent with ZFC that $\rr<\rr_{fi}$.
\end{cor}

\begin{pf}
  Cohen's original model for the negation of the continuum hypothesis
  has $\non M=\aleph_1$ and $\cov M=\ger c$. It follows, by
  Theorems~\ref{belownonM:thm} and \ref{abovecovM:thm} that this model
  satisfies
\[
\rr\leq\non M=\aleph_1<\ger c=\cov M\leq\rr_{fi}.  \qedhere
\]
\end{pf}

\section{Bounding and Dominating, More Padding} \label{pad2}

In this section, we extend the method of padding with zeros, used in
the proof of Theorem~\ref{uncbl:thm}, to obtain stronger lower bounds
for rearrangement numbers.  Recall that the key idea in the padding
method was to spread out the nonzero terms in a series so far that the
permutations under consideration do not change their relative order
(up to finitely many exceptions).  The following definition introduces
a cardinal characteristic intended to capture this idea.

\begin{df}      \label{j:df}
A set $A\subseteq\bbb N$ is \emph{preserved} by a permutation $p$ of
\bbb N if $p$ does not change the relative order of members of $A$
except for finitely many elements.  That is, for all but finitely many
elements of $A$, we have $x<y\iff p(x)<p(y)$.  If $A$ is not
preserved by $p$ we say that it is \emph{jumbled} by $p$.  A
\emph{jumbling family} is a family of permutations such that every
infinite $A\subseteq\bbb N$ is jumbled by at least one member of the
family.  The \emph{jumbling number}, \ger j, is the smallest
cardinality of a jumbling family.
\end{df}

The proof of Theorem~\ref{uncbl:thm} shows that \ger j is
uncountable.  We shall see later that $\ger j=\ger b$, but first we
check that the concept of jumbling provides a lower bound for the
rearrangement numbers.

\begin{thm}     \label{abovej:thm}
$\ger j\leq\rr$.
\end{thm}

\begin{pf}
Let $C$ be any family of fewer than \ger j permutations; we shall find
a conditionally convergent series whose sum is unchanged under all the
permutations in $C$.  Let $\sum_nb_n$ be any conditionally convergent
series, and, since $|C|<\ger j$, let $A\subseteq\bbb N$ be an infinite
set preserved by the inverses of all the permutations in $C$.  Let
$\sum_na_n$ be the series obtained by putting the $b_k$'s at positions
in $A$, in order, and filling the remaining positions with zeros.
That is,
\[
a_n=
\begin{cases}
  b_k&\text{if }n\text{ is the  }k\th\text{ element of }A,\\
0&\text{if }n\notin A.
\end{cases}
\]
If $p\in C$ then, since $A$ is preserved by $p^{-1}$, the orders of the
nonzero terms of $\sum_na_n$ and of its rearrangement $\sum_na_{p(n)}$
are the same, with at most finitely many exceptions.  Therefore, the
sums agree.  We have a conditionally convergent series, $\sum_na_n$,
whose sum is unchanged when it is rearranged by any of the
permutations in $C$.
\end{pf}

To connect this theorem with familiar cardinal characteristics, we
show next that the jumbling number is the same as the unbounding
number.

\begin{thm}     \label{aboveb:thm}
  $\ger j=\ger b$.  Consequently, $\ger b\leq\rr$.
\end{thm}

\begin{pf}
  We need only prove $\ger j=\ger b$, because the ``consequently''
  part of the theorem then follows immediately via
  Theorem~\ref{abovej:thm}.

  We begin by proving $\ger b\leq\ger j$, i.e., a family $C$ of fewer
  than \ger b permutations cannot be a jumbling family; it must
  preserve some infinite set.  Consider any such family $C$ and
  associate to each permutation $p\in C$ a function
  $f_p:\bbb N\to\bbb N$ with the property that, for every
  $n\in\bbb N$, we have $n<f_p(n)$ and
\[
(\forall x\leq n)(\forall y\geq f_p(n))\,p(x)<p(y).
\]
Notice that we are only requiring $f_p(n)$ to be larger than finitely
many numbers, namely $n$ and all the numbers $p^{-1}(z)$ for
$z\leq\max(p[n+1])$. So there is no difficulty obtaining such a
function $f_p$.

Because $|C|<\ger b$, there is a function $g:\bbb N\to\bbb N$ such
that $f_p\leq^*g$ for all $p\in C$.  Increasing the values of $g$, we
can arrange that $g$ is strictly increasing.  Now define an infinite
set $A=\{a_0<a_1<\dots\}$ of natural numbers inductively, starting
with an arbitrary $a_0$ (say 0) and ensuring at every step that
$a_{n+1}\geq g(a_n)$.  We claim that all the permutations $p\in C$
preserve $A$.

To see this, consider any $p\in C$ and, since $f_p\leq^*g$, fix $k$
such that $f_p(n)\leq g(n)$ for all $n\geq k$.  For elements $a_r<a_s$
of $A$ that are larger than $k$, we have, since $r+1\leq s$,
\[
f_p(a_r)\leq g(a_r)\leq a_{r+1}\leq a_s,
\]
and therefore, by our choice of $f_p$, $p(a_r)<p(a_s)$.  This shows
that every $p\in C$ preserves $A$, and so it completes the proof that
$\ger b\leq\ger j$.  (Note that this inequality suffices to give the
``consequently'' part of the theorem.)

It remains to show that $\ger j\leq\ger b$.  For this purpose, it is
convenient to invoke an alternative characterization of \ger b,
essentially due to Solomon~\cite{sol}; for the version used here, see
\cite[Theorem~2.10]{hdbk}.  An \emph{interval partition} is a
partition of \bbb N into (infinitely many) finite intervals
$I_n=[i_n,i_{n+1})$, where $0=i_0<i_1<\dots$. A second interval
partition $\{J_n:n\in\bbb N\}$ is said to \emph{dominate} the interval
partition $\{I_n:n\in\bbb N\}$ if, for all but finitely many $k$, the
interval $J_k$ includes some $I_n$ as a subset.  Then \ger b is the
smallest cardinality of a family of interval partitions such that no
single interval partition dominates them all.

Fix an undominated family \scr F of \ger b interval partitions.  To
each of the partitions $I=\{I_n:n\in\bbb N\}\in\scr F$, associate the
permutation $p_I$ that flips each of the intervals $I_n$ upside down.
That is, if $x\in I_n=[i_n,i_{n+1})$ then $p_I(x)=i_n+i_{n+1}-x-1$.
Let $C=\{p_I:I\in\scr F\}$. So $C$ has cardinality \ger b, and we
shall complete the proof by showing that $C$ is a jumbling family.

Consider any infinite $A=\{a_0<a_1<\dots\}\subseteq\bbb N$, and assume
without loss of generality that $a_0=0$. Form an interval partition
$J=\{J_k:k\in\bbb N\}$ by setting $J_k=[a_{3k},a_{3k+3})$.  By our
choice of \scr F, it contains an interval partition $I=\{I_n:n\in\bbb
N\}$ that is not dominated by $J$.  That is, for infinitely many
values of $k$, there is no interval $I_n$ included in $J_k$.

Temporarily fix one such $k$, and let $I_n$ be the interval of
$I$ that contains $a_{3k+1}$. If $I_n$ contained neither $a_{3k}$ nor
$a_{3k+2}$, then $I_n$ would be included in $J_k$, contrary to our
choice of $k$.  Therefore, $I_n$ must contain at least one element of
$A$ in addition to $a_{3k+1}$.  We thus have two elements of $A$ whose
order is reversed by $p_I$, since $p_I$ flips $I_n$ upside down.

Now un-fix $k$. The preceding paragraph applies to infinitely many
values of $k$, so we have infinitely many pairs of elements of $A$
whose order is reversed by $p_I$.  So $p_I$ jumbles $A$.  Since $A$
was arbitrary, $C$ is a jumbling family.
\end{pf}

Since \rr\ is the smallest of the rearrangement numbers, any lower
bound for it, such as \ger b from Theorem~\ref{aboveb:thm},
automatically applies to the variants $\rr_{fi}$, $\rrf$, and \rri\ as
well.  These  three variants, however, admit the following stronger
lower bound.

\begin{thm}     \label{aboved:thm}
$\ger d\leq\rr_{fi}$. Consequently, $\ger d\leq\rrf$ and $\ger
d\leq\rri$.
\end{thm}

\begin{pf}
  We need only prove that $\ger d\leq\rr_{fi}$, since the
  ``consequently'' part of the theorem then follows by virtue of
  the trivial inequalities pointed out in Section~\ref{def}.

Consider any family $C$ of fewer than \ger d permutations of \bbb N.
We must find a conditionally convergent series $\sum_na_n$ such that
none of the permutations in $C$ make this series converge to a
different finite sum or diverge to $+\infty$ or $-\infty$.

Begin by associating to each $p\in C$ a function $f_p:\bbb N\to\bbb N$
as in the proof of Theorem~\ref{aboveb:thm} except that we use
$p^{-1}$ instead of $p$.  That is, for every $n\in\bbb N$, we have
$n<f_p(n)$ and
\[
(\forall x\leq n)(\forall y\geq f_p(n))\,p^{-1}(x)<p^{-1}(y).
\]
Because $|C|<\ger d$, there exists $g:\bbb N\to\bbb N$ that is not
eventually dominated by any of these functions $f_p$ for $p\in C$.  We
can arrange, by increasing its values if necessary, that $g$ is
strictly increasing and that $g(n)>n$ for all $n$.  By iterating $g$,
we obtain a strictly increasing sequence
\[
0<g(0)<g(g(0))<\dots<g^k(0)<g^{k+1}(0)<\dots.
\]
This sequence will be used to apply the padding-with-zeros technique as
follows.

Fix a conditionally convergent series $\sum_nb_n$, and define the
padded version $\sum_na_n$ by
\[
a_n=
\begin{cases}
  b_k & \text{if }n=g^k(0)\\
0 & \text{if }n\text{ is not of the form }g^l(0)\text{ for any }l.
\end{cases}
\]
We shall show that no $p\in C$ can have $\sum_na_{p(n)}$ converging to
a finite sum other than $\sum_na_n$, nor can that rearranged sum
diverge to $+\infty$ or to $-\infty$.

Consider any $p\in C$. By our choice of $g$, there are infinitely many
numbers $x\in\bbb N$ such that $f_p(x)<g(x)$.  Temporarily concentrate
on one such $x$. Let $k$ be the smallest integer such  that
$x<g^k(0)$.  So we have both $g^{k-1}(0)\leq x$ and $g(x)<g^{k+1}(0)$
(the former because $k$ is smallest and the latter because $g$ is
increasing).  Therefore,
\[
g^{k-1}(0)\leq x<f_p(x)<g(x)<g^{k+1}(0),
\]
where the second inequality comes from our choice of $f_p$ and the
third from our choice of $x$.  In view of the definition of $f_p$, we
see that all the numbers
$p^{-1}(0),p^{-1}(g(0)),\dots,p^{-1}(g^{k-1}(0))$ are smaller than all
the numbers $p^{-1}(g^{k+1}(0)),p^{-1}(g^{k+2}(0)),\dots$. Notice
that, for any $r$, $p^{-1}(g^r(0))$ is the location where $b_r$
appears in the series $\sum_na_{p(n)}$.  so in this series, all of
$b_0,b_1,\dots,b_{k-1}$ occur before all of $b_{k+1},b_{k+2},\dots$.
(Nothing is said here about $b_k$; it could occur out of order
anywhere.)  Therefore, this rearranged series has a partial sum that
differs from $\sum_{n=0}^kb_n$ by at most $|b_k|$.

This discussion was based on a particular $x$ where $f_p(x)<g(x)$.
But there are infinitely many such $x$'s and infinitely many $k$'s
associated to them as above.  For each of these $k$'s, we have seen
that $\sum_na_{p(n)}$ has a partial sum differing from
$\sum_{n=0}^kb_n$ by at most $|b_k|$.  But, as $k$ tends to infinity,
$\sum_{n=0}^kb_n$ tends to the infinite sum $\sum_nb_n=\sum_na_n$, and
$|b_k|$ tends to zero.  Thus, in the sequence of partial sums of
$\sum_na_{p(n)}$, there is an infinite subsequence tending to
$\sum_na_n$, which means that the whole sequence of partial sums
cannot tend to a different value or to $\pm\infty$.
\end{pf}

Note that the argument at the end of this proof does not contradict
the possibility that $\sum_na_{p(n)}$ might diverge by oscillation,
with its partial sums $\sum_{n=0}^ka_{p(n)}$ coming close to
$\sum_na_n$ at the infinitely many $k$'s under consideration, but
wandering around for other values of $k$.  This is why the theorem
gives a lower bound for $\rr_{fi}$ but not for \rr.  In fact, it is
not provable in ZFC that $\ger d\leq\rr$.  This is a consequence of
Theorem~\ref{belownonM:thm} and the fact that $\non M<\ger d$ is known
to be consistent with ZFC. (The basic Cohen model satisfies
$\non M=\aleph_1<\ger d=\ger c$.)

Note also that Theorem~\ref{aboved:thm} supersedes
Theorem~\ref{abovecovM:thm} because it is provable that
$\cov M\leq\ger d$, and it is consistent that this inequality is
strict. (For example, strict inequality holds in the Laver and Miller
models.)

\section{Measure, Random Signs}         \label{signs}

In this section, we relate the rearrangement numbers to the covering
number for measure.  We shall need a result of Rademacher~\cite{rade},
stated as a lemma below, about infinite series with randomly chosen
signs. This lemma can be viewed as a special case of Kolmogorov's famous three-series theorem, which was proved later.

\begin{la}[Rademacher]
Let $(c_n:n\in\bbb N)$ be any sequence of real numbers.  Let
$A\subseteq2^{\bbb N}$ be the set of those sequences $s$ of zeros and
ones for which $\sum_n(-1)^{s(n)}c_n$ converges.  Then the Lebesgue
measure of $A$ is 1 if $\sum_n{c_n}^2$ converges and 0 otherwise.
\end{la}

In other words, if we attach signs randomly to the terms of the series
$\sum_nc_n$, the result will converge almost surely if $\sum_n{c_n}^2$
converges, and it will diverge almost surely otherwise.

\begin{thm}     \label{abovecovN:thm}
$\cov N\leq\rr$.
\end{thm}

\begin{pf}
Consider first a single permutation $p$ of \bbb N.  Since
$\sum_n1/p(n)^2$ is a rearrangement of the convergent series
$\sum_n1/n^2$, it converges, and therefore, by the lemma, the set
\[
A_p=\{s\in2^{\bbb N}:\sum_n(-1)^{s(n)}/p(n)\text{ diverges}\}
\]
has measure zero.  Therefore, so does
\[
B_p=\{s\in2^{\bbb N}:\sum_n(-1)^{s(p(n))}/p(n)\text{ diverges}\},
\]
since it is just the pre-image of $A_p$ under the measure-preserving
bijection $s\mapsto s\circ p:2^{\bbb N}\to2^{\bbb N}$.

Now consider any family $C$ of fewer than \cov N permutations of \bbb
N.  By definition of \cov N, the associated measure-zero sets $B_p$
for $p\in C$ cannot cover $2^{\bbb N}$, so there is some
  $s\in2^{\bbb N}$ such that $\sum_n(-1)^{s(p(n))}/p(n)$ converges for
  all $p\in C$.  That is, the rearrangements of $\sum_n(-1)^{s(n)}/n$
  by permutations from $C$ will not diverge to $\pm\infty$ or by
  oscillation.

This proves that $\cov N\leq\rr_{io}$, and it remains only to recall
from Theorem~\ref{mix-o:thm} that $\rr_{io}=\rr$.
\end{pf}

\begin{rmk}
  Readers familiar with Cicho\'n's diagram of cardinal charactersitics
  (see for example \cite[Section~5]{hdbk}) will notice that the bounds we have proved
  for \rr, namely $\max\{\cov N,\ger b\} \leq\rr\leq\non M$, given by
  Theorems~\ref{belownonM:thm}, \ref{aboveb:thm}, and \ref{abovecovN:thm}, sandwich \rr\
  between adjacent characteristics in that diagram.
\end{rmk}

There is no provable inequality in either direction between \cov N and
\ger b.  Specifically, $\cov N<\ger b$ in the Laver model and
$\ger b<\cov N$ in the random real model.  Thus, the lower bounds for
\rr\ in Theorems~\ref{aboveb:thm} and \ref{abovecovN:thm} are
independent, and each of them can consistently be strict.

Furthermore, since the random model has $\ger d<\cov N$, we find that
the lower bounds \ger d and \cov N for $\rr_{fi}$ are independent and
can each consistently be strict.  We summarize these consistency
observations for future reference.

\begin{cor}     \label{con-strict:cor}
None of the inequalities
\begin{ls}
  \item $\cov N\leq\rr$,
\item $\ger b\leq\rr$,
\item $\cov N\leq\rr_{fi}$,
\item $\ger d\leq\rr_{fi}$
\end{ls}
is provably reversible.  That is, in each case, strict inequality is
consistent with ZFC.
\end{cor}

\section{Forcing New Finite Limits}    \label{force1}

In this section and the next, we show that it is consistent with ZFC
that all the rearrangement numbers are strictly smaller than the
cardinality of the continuum. In view of the order relationships
between the various rearrangement numbers (see Figure~\ref{fig2}), it
suffices to prove the consistency of the two inequalities $\rrf<\ger
c$ and $\rri<\ger c$.  In this section, we construct a model for
$\rrf<\ger c$.  In Section~\ref{force2}, we shall construct a model
for $\rri<\ger c$, and we shall point out that the two constructions
can be combined to produce a model that satisfies both inequalities
simultaneously.  In these two sections, we assume a familiarity with the
forcing method of building models of ZFC, including the technique of
finite-support iterated forcing.

Both models will be obtained by starting with a model in which
$\ger c>\aleph_1$ (we can take $\ger c$ as large as we wish), and then
performing an $\omega_1$-step, finite-support iteration of forcings
that satisfy the countable chain condition.  Because of the chain
condition, cardinals in the final model will be the same as in the
ground model, so the final model will have $\ger c>\aleph_1$.  At each
step of the iteration, we shall adjoin a permutation of \bbb N that
disrupts the convergence of all conditionally convergent series in the
model produced by the previous steps; in the present section, the
disruption consists of producing a new, finite sum for the rearranged
series, and in Section~\ref{force2} it consists of making the
rearranged series diverge to $+\infty$ or to $-\infty$.  Thanks to the
countable chain condition, every conditionally convergent series in
the final model is already in the intermediate model after some
countably many steps, and so its convergence will be disrupted by the
permutation added at the next step.  Thus, the $\aleph_1$ permutations
that we added, one per step of the iteration, suffice to disrupt the
convergence of all conditionally convergent series in the final
model. That is, the final model will satisfy $\rrf<\ger c$ in the
construction from the
present section, $\rri<\ger c$ in the construction from
Section~\ref{force2}, and both inequalities in a construction that
interleaves the two iterations.

We turn now to the construction of a forcing notion that satisfies the
countable chain condition and adds a permutation of \bbb N that
rearranges all conditionally convergent series in the ground model to
have new, finite sums.  In fact, the new sums will be very new, in that
they are outside the ground model.  For this proof, we shall need a
classical result of L\'evy~\cite{levy} and Steinitz~\cite{steinitz}
and the following associated definitions.

\begin{df}
  Let $d$ be a natural number, and let $\bar a=\sq{a^i:i<d}$ be a
  $d$-tuple of infinite series of real numbers. Define $K(\bar a)$ to
  be the set of $d$-tuples of real numbers \sq{s_i:i<d} for which the
  series $\sum_{i<d}s_ia^i$ converges absolutely.  Define $R(\bar a)$
  to be the orthogonal complement of $K(\bar a)$ in $\bbb R^d$.  The
  $d$-tuple of series \sq{a^i:i<d} is said to be \emph{independent} if
  $K(\bar a)=\{0\}$.  An arbitrary set $I$ of infinite series of real
  numbers is said to be \emph{independent} if every finite tuple of
  distinct elements of $I$ is independent.
\end{df}

Notice that independence as defined here is the ordinary notion of
linear independence applied to the quotient of the vector space of all
infinite series of real numbers modulo the subspace of absolutely
convergent series.

The L\'evy-Steinitz Theorem extends the Riemann Rearrangement Theorem
to the context of infinite series of vectors in a finite-dimensional
space $\bbb R^d$.

\begin{thm}[L\'evy \cite{levy}, Steinitz
  \cite{steinitz}] \label{ls:thm}
  If $\bar a=\sq{a^i:i<d}$ is a finite tuple of convergent series of
  real numbers, then the set of $d$-tuples \sq{\sum_na^i_{p(n)}:i<d}
  obtainable by rearrangements $p$ coincides with the set of vector
  sums $\sq{\sum_na^i_n:i<d}+\bar x$ with $\bar x\in R(\bar a)$.
\end{thm}

In other words, the alterations of the sum \sq{\sum_na_n^i:i<d}
obtainable by permuting the summands are the same as the alterations
that simply add an arbitrary vector from $R(\bar a)$.

More modern sources than \cite{levy} and \cite{steinitz} for
information about the L\'evy-Steinitz Theorem include
\cite{bd,kadets,rosen}.

We shall make use of this theorem via the following corollary.

\begin{cor}     \label{ls:cor}
  Let $\bar a=\sq{a^i:i<d}$ be an independent $d$-tuple of convergent
  series of real numbers. Let $\bar v$ be an arbitrary vector in
  $\bbb R^d$.  Let $f:n\to\bbb N$ be an injective function from some
  natural number $n$ into \bbb N.  Then there is a permutation $p$ of
  \bbb N, extending $f$, such that $\sq{\sum_na^i_{p(n)}:i<d}=\bar v$.
\end{cor}

\begin{pf}
Since $\bar a$ is independent, we have $R(\bar a)=\bbb R^d$, so $\bar
v$ is, by the L\'evy-Steinitz theorem, obtainable as the rearranged sum
\sq{\sum_na^i_{p(n)}:i<d} for some permutation $p$ of \bbb N.  To make
$p$ extend $f$, it suffices to alter $p$ at only finitely many places,
and this will not affect the sum of the rearrangement.
\end{pf}

We shall also need the following result of Steinitz, which is used in
one of the proofs of the L\'evy-Steinitz Theorem.

\begin{thm}[Polygonal Confinement,
  \cite{steinitz}]     \label{poly:thm}
For each natural number $d$, there exists a constant $C_d$ such that,
if $\bar v_0,\bar v_1,\dots,\bar v_{n-1}$ are any $n$ vectors in $\bbb
R^d$, each of length $\nm{\bar v_m}\leq1$, with sum zero,
$\sum_{m<n}\bar v_m=0$, then there is a permutation $p$ of
$n\setminus\{0\}$ such that
\[
\nmb{\bar v_0+\sum_{m\in k\setminus\{0\}}\bar v_{p(m)}}\leq C_d
\]
for all $k\leq n$.
\end{thm}

In other words, given a closed polygonal path, starting and ending at
the origin in $\bbb R^d$, whose sides have lengths $\leq1$, one can
reorder the sides (except for the first) so that the entire polygon
stays within $C_d$ of the origin.  The essential point is that $C_d$
depends only on the dimension, not on the number of  steps $\bar
v_m$ in the path.

Fix, once and for all, a nondecreasing sequence of constants $C_d$
satisfying the conclusion of the theorem.

The following corollary is a slight variant of the theorem and will be
more convenient for our application.

\begin{cor}     \label{poly:cor}
Let $\bar v_0,\bar v_1,\dots,\bar v_{n-1}$ be $n$ vectors in $\bbb
R^d$, let $\bar b$ be their sum, and let $\rho$ be a positive real
number with all $\nm{\bar v_i}\leq\rho$ and $\nm{\bar b}\leq\rho$.
Then there is a permutation $p$ of $n\setminus\{0\}$ such that
\[
\nmb{\bar v_0+\sum_{i\in k\setminus\{0\}}\bar v_{p(i)}}\leq\rho
C_d+\nm{\bar b}
\]
for all $k\leq n$.
\end{cor}

\begin{pf}
  Dividing all the vectors $\bar v_i$ and their sum $\bar b$ by
  $\rho$, we may assume without loss of generality that $\rho=1$.
  Then, since the sequence
  $\bar v_0,\bar v_1,\dots,\bar v_{n-1}, -\bar b$ has sum zero and
  consists of vectors of length at most 1, we can apply the Polygonal
  Confinement Theorem to obtain a permutation of this sequence in
  which $\bar v_0$ is still first, and all partial sums (starting at
  the beginning of the sequence) have length at most $C_d$.  If
  $-\bar b$ were still at the end of the sequence, after this
  permutation, then these partial sums would be the partial sums in
  the statement of the corollary, and we would be done (even without
  the $\nm{\bar b}$ term on the right side of the inequality).  If
  $-\bar b$ is not at the end after the permutation, then some of the
  partial sums that we know to be shorter than $C_d$ would differ from
  the partial sums in the statement of the corollary; the former would
  include $-\bar b$ while the latter would not.  But that difference
  affects the lengths of these partial sums by at most \nm{\bar b},
  by the triangle inequality.  So we get the inequality claimed in the
  corollary.
\end{pf}

We are now ready to define the partial order \bbb P that will add a
permutation making all conditionally convergent series from the ground
model converge to new (finite) sums not in the ground model.

\begin{conv}
  Throughout this section, $I$ is an independent set of convergent
series of real numbers.
\end{conv}

Independence is not needed for the
definition of \bbb P, but it is involved in the subsequent lemmas
showing that \bbb P has the desired properties.  Note that
independence implies that all the series in $I$ are conditionally
convergent.

\begin{df}
$\bbb P_I$ is the partially ordered set whose elements are triples
$(f,A,\eps)$ such that:
\begin{ls}
\item$f$ is an injective function from some $n\in\bbb N$ into \bbb N.
\item $A$ is a finite nonempty subset of $I$.
\item $\eps$ is a positive rational number.
\item If \sq{a^i:i<d} is an enumeration of $A$, then
  \[
\nm{\sq{a^i_m:i<d}} <\eps/C_d
\]
for all $m\in\bbb N\setminus\ran f$.
\end{ls}
The order on $\bbb P_I$ is defined by setting
$(g,B,\delta)\leq(f,A,\eps)$
when:
\begin{ls}
\item $g$ extends $f$.
\item $B$ is a superset of $A$.
\item If \sq{a^i:i<d} enumerates $A$, then
\begin{itemize}
\item[$\circ$]
for all $m \in \dom g+1$,
$$   
\nmb{\sum_{k\in m\setminus\dom f}\sq{a^i_{g(k)}:i<d}}<\eps.
$$
\item[$\circ$]
\vspace{1mm}
$\displaystyle \quad \ 
\delta+\nmb{\sum_{k\in\dom g\setminus\dom f}\sq{a^i_{g(k)}:i<d}}\leq\eps.
$
\end{itemize}
\end{ls}
\end{df}

\begin{rmk}
This remark is an attempt to aid the reader's intuition about this
notion of forcing; it can be skipped by those readers who are willing
to simply work with the formal definition of $\bbb P_I$.

In any condition $(f,A,\eps)$, the first component $f:n\to\bbb N$ is
intended to be an initial segment of the generic permutation $\pi$ added
by the forcing.  Thus, the first $n$ terms of a rearranged series
$\sum_nt_{\pi(n)}$ will be the terms of the original series in the
  locations specified by $f$, namely $t_{f(0)},\dots,t_{f(n-1)}$.

  The second component, $A$, specifies finitely many series
  $\sum_ma^i_m$ in $I$ over which our condition wants to exercise some
  control.  The last clause in the definition of conditions says that,
  except for those terms whose position in the $\pi$-rearrangement has
  already been specified by $f$, the remaining terms in the series in
  $A$ are small compared to $\eps$.  In fact they are very small in
  two senses.  First, the inequality applies to these terms not just
  individually but ``jointly'' across all elements of $A$. That is, it
  does not just bound the individual terms $a^i_m$  but the
  $d$-component vectors \sq{a^i_m:i<d}.  Second, the bound is not
  merely $\eps$ but $\eps/C_d$, where $C_d$ is the constant from the
  L\'evy-Steinitz Theorem.  The point of this is that it provides, via
  Corollary~\ref{poly:cor}, a bound for sums of these vectors if we are
  willing to suitably rearrange them.

  The use of an enumeration \sq{a^i:i<d} in the last requirement for
  conditions, and also later in the definition of the ordering and
  elsewhere, is unimportant in the sense that, if the statements are
  true for one enumeration of $A$, then they are also true for all
  other enumerations.  The only reason enumerations are involved at
  all is to have an ordering of the components in vectors like
  \sq{a^i_m:i<d}.  If we stretched the meaning of ``vector'' to allow
  the components to be indexed by finite sets other than natural
  numbers, then no such enumeration would be needed; $A$ itself could
  serve as the index set.

  In the definition of the ordering, the first two clauses are
  standard; a stronger condition tells us more about the generic
  permutation $\pi$ (i.e., it specifies a longer initial segment of
  $\pi$), and it tries to control more of the series in $I$.  The last
  two clauses are more subtle, but it helps to notice first that they
  refer only to the series in $A$, the ones that the weaker condition
  $(f,A,\eps)$ wants to control. $B$ is not mentioned in these
  clauses.  Furthermore, these clauses are about the vectors
  \sq{a^i_q:i<d} associated to locations $q$ in
  $\dom g\setminus\dom f$, i.e., locations for which $f$ did not say
  where they will go in the $\pi$-rearrangement but $g$ did.  If we
  think of these vectors as listed in a sequence, in the order
  assigned to them by $g$, then all initial segments of this sequence
  are required to have small sums, i.e., shorter than $\eps$; so,
  intuitively, $g$ arranged these vectors in an intelligent order, as
  suggested by polygonal confinement.  And furthermore, the amount by
  which the final sum of all these vectors is shorter than $\eps$ is
  an upper bound for the third component $\delta$ in the stronger
  condition.  The point of that is that further extensions will be
  subject to bounds given by this $\delta$ and that will prevent them
  from combining with the extension given by $g$ to achieve sums
  greater than $\eps$.
\end{rmk}

It is not difficult to check that the definition of the ordering of
$\bbb P_I$ is legitimate; it is, in particular, transitive.  (For some
intuition behind transitivity, see the last sentence of the preceding
remark.)  It is also not difficult to see that $\bbb P_I$ is
nonempty. In fact, for any nonempty finite subset $A$ of $I$, there is
an $\eps$ such that $(\emp,A,\eps)$ is a condition.  To prove it, use
the fact that all the series in $A$ converge, so their terms are
bounded, and then just choose $\eps$ large enough.

We next prove several lemmas establishing density properties of $\bbb
P_I$.  All of them depend on our convention that $I$ is an independent
family of convergent series.  The first lemma lets us extend the
initial segment $f$ of the generic permutation and tighten the
constraint $\eps$.

\begin{la}      \label{dense1}
  For any condition $(f,A,\eps)$ and any positive integer $n$, there
  is an extension $(g,A,\delta)\leq(f,A,\eps)$ (with the same second
  component $A$) with the following
  properties:
  \begin{ls}
    \item $n\subseteq\dom g\cap\ran g$.
\item $\delta<1/n$.
\item If \sq{a^i:i<d} is an enumeration of $A$ then
\[
\nm{\sq{a^i_m:i<d}}<\frac\delta{2C_{2d}}
\]
for all $m\in\bbb N\setminus\ran g$.
  \end{ls}
\end{la}

\begin{pf}
  Let $(f,A,\eps)$ and $n$ be given, let $r=\dom f$, and fix an
  enumeration \sq{a^i:i<d} of $A$.  By Corollary~\ref{ls:cor}, there
  is a permutation $p$ of \bbb N, extending $f$, such that
\[
\sum_{n=0}^\infty\sq{a^i_{p(n)}:i<d}=\sum_{n\in\dom
  f}\sq{a^i_{f(n)}:i<d},
\]
or in other words,
\[
\sum_{m\geq r}\sq{a^i_{p(m)}:i<d}=0.
\]
Recall that, by definition of conditions, we have strict inequalities
$\nm{\sq{a^i_m:i<d}}<\eps/C_d$ for all $m\notin\ran f$.  Furthermore,
the norms on the left side of these inequalities tend to zero as $m$
increases, because the series in $A$ are convergent.  So we can fix a
positive number $\eta<\eps$ such that $\nm{\sq{a^i_m:i<d}}<\eta/C_d$
for all $m\notin\ran f$.  Fix a positive rational number $\delta$
smaller than both $1/n$ and$(\eps-\eta)/2$.

For any sufficiently large natural number $n_*$, we have all of the following:
\begin{ls}
  \item $n_*\geq n$.
\item $n\subseteq p[n_*]$.
\item $\nmb{\sum_{m\in n_*\setminus r}\sq{a^i_{p(m)}:i<d}}<\delta$.
\item $\nm{\sq{a^i_m:i<d}}<\delta/2C_{2d}$ for each $m\in\bbb
  N\setminus n_*$.
\end{ls}
The first and second of these assertions are clear, and the fourth
follows from the fact that all the series in $A$ converge. To see the
third, note that the sum there is a partial sum of the series
$\sum_{m\geq r}\sq{a^i_{p(m)}:i<d}$ whose sum is zero by our
choice of $p$.

Choose $n_*$ large enough so that all these statements are true, and
then use Corollary~\ref{poly:cor} to produce an injection from
$n_*\setminus r$ to \bbb N which, when combined with $f:r\to\bbb N$,
produces an injection\footnote{Corollary~\ref{poly:cor} gives us
  $g(r)=r$ here. This information, though used in the formulas that
  follow, is not essential for the proof here or in
  Lemma~\ref{dense2}.}  $g:n_*\to\bbb N$, extending $f$, with the same
range as $p\restr n_*$, such that, for all $m\leq n_*$,
\[
\nmb{\sq{a^i_r:i<d}+\sum_{k\in m\setminus
    (r+1)}\sq{a^i_{g(k)}:i<d}}\leq (\eta/C_d)C_d+\delta<\eps-\delta.
\]
Then $(g,A,\delta)$ is as required in the lemma.
\end{pf}

The next lemma allows us to enlarge the set $A$ of controlled series.

\begin{la}      \label{dense2}
  For each condition $(f,A,\eps)\in\bbb P_I$ and each $b\in I$, there
  is an extension $(g,B,\delta)\leq(f,A,\eps)$ with $b\in B$.
\end{la}

\begin{pf}
  Assume $b\notin A$, as otherwise there is nothing to prove.
  Enumerate $A\cup\{b\}$ as \sq{a^i:i\leq d} with $b$ as the last
  element in the enumeration, $b=a^d$.  As in the proof of the
  preceding lemma, use Corollary~\ref{ls:cor} to extend $f$ to a
  permutation $p$ of \bbb N such that $\sum_{m\geq
    r}\sq{a^i_{p(m)}:i<d}=0$, where, as before, $r$ is the domain of
  $f$.  Continuing as in the previous proof, choose $\eta$, $\delta$,
  and $n_*$ as there except that the fourth condition on $n_*$ is
  strengthened to include $b$ with the other $a^i$'s and weakened by
  using $C_{d+1}$ in place of $C_{2d}$, i.e.,
  \begin{ls}
  \item $\nm{\sq{a^i_m:i\leq d}}<\delta/C_{d+1}$ for each
    $m\in\bbb N\setminus n_*$.
  \end{ls}
The strengthening is easy to obtain because $b$ as well as the other
$a^i$'s are convergent series.

Finally, still proceeding as in the previous proof but with $b$
included, use Corollary~\ref{poly:cor} to extend $f$ to an injection
$g:n_*\to\bbb N$ such that, for all $m\leq n_*$,
\[
\nmb{\sq{a^i_r:i<d}+\sum_{k\in m\setminus
    (r+1)}\sq{a^i_{g(k)}:i<d}}\leq
(\eta/C_{d+1})C_{d+1}+\delta<\eps-\delta.
\]
Then $(g,A\cup\{b\},\delta)$ is as required in the lemma.
\end{pf}

The preceding two lemmas provide the following important information
about the generic object added by forcing with $\bbb P_I$.

\begin{cor}     \label{perm-conv:cor}
If $G\subseteq\bbb P_I$ is a $V$-generic filter and we define
\[
\pi=\bigcup\{f:(f,A,\eps)\in G\},
\]
then $\pi$ is a permutation of \bbb N and, for every series $a\in I$,
the rearrangement $\sum_na_{\pi(n)}$ converges.
\end{cor}

\begin{pf}
Since all the first components $f$ of conditions in $G$ are
single-valued, injective, and pairwise compatible, $\pi$ is a partial
function from \bbb N to \bbb N. That it is total and surjective, and
thus a permutation of \bbb N, follows from genericity and the clause
$n\subseteq\dom g\cap\ran g$ in Lemma~\ref{dense1}.

For any series $a\in A$, genericity and Lemma~\ref{dense2} provide a
condition $(f,A,\eps)\in  G$ with $a\in A$; by Lemma~\ref{dense1}, we
can further arrange that $\eps$ here is as small as we want.  Then, by
definition of the ordering of $\bbb P_I$, extensions of $(f,A,\eps)$
cannot produce large variations in the partial sums of
$\sum_{n=r}^\infty a_{\pi(n)}$, where $r=\dom f$.  Since any partial sum of this generic
rearrangement is obtainable from a condition in $G$, it is also
obtainable from an extension of $(f,A,\eps)$.  So these partial sums
cannot oscillate by more than $\eps$.  Since $\eps$ can be taken to
be as small as we want, it follows that $\sum_na_{\pi(n)}$ converges.
\end{pf}

The next (and last) of the density lemmas serves to ensure that the
sum of the series rearranged by $\pi$ is not in the ground model.

\begin{rmk}     \label{indep:rmk}
Any effort to impose a particular behavior (in the present situation,
the behavior of convergence to new values) on arbitrary conditionally
convergent series must confront the fact that two or more series might
be related in such a way that their behavior under rearrangements is
correlated, possibly in undesirable ways.  Until now, the present
argument has avoided this issue by dealing with an independent set $I$
of series.  Utimately, though, it will have to deal with arbitrary
series in the ground model.  The following lemma is a key step in this
direction, dealing with linear combinations of series from $I$.
Later, by taking $I$ to be a \emph{maximal} independent set, we shall
use this lemma to deal with all series in the ground model.
\end{rmk}

\begin{la}      \label{dense3}
Let $(f,A,\eps)$ be a condition in $\bbb P_I$, let \sq{a^i:i<d} be an
enumeration of $A$, let \sq{s_i:i<d} be a $d$-tuple of nonzero real
numbers, and let $r$ be any real number. Then there exists an
extension $(g,A,\delta)\leq(f,A,\eps)$ such that
\[
\big|r-\sum_{i<d}\sum_{n\in\dom
    g}s_ia^i_{g(n)}\big|>\delta\sum_{i<d}|s_i|.
\]
\end{la}

\begin{pf}
As a preliminary step, we extend the given condition, if necessary, to
obtain
\[
r\neq\sum_{i<d}\sum_{n\in\dom f}s_ia^i_{f(n)}.
\]
If the desired inequality does not already hold, then we proceed as
follows. Since $I$ is independent and the $s_i$ are nonzero, the
series $\sum_n\left(\sum_{i<d}s_ia^i_n\right)$ is conditionally
convergent.  So there are arbitrarily large $m$ with
$\sum_{i<d}s_ia^i_m\neq0$. As in previous lemmas, choose $\eta<\eps$
such that $\nm{\sq{a^i_m:i<d}}<\eta/C_d$ for all $m\notin\ran f$, and
let $\delta$ be a positive rational number smaller than
$(\eps-\eta)/2$.  Then find an $m\notin\ran f$ such that both
$\sum_{i<d}s_ia^i_m\neq0$ and $\nm{\sq{a^i_m:i<d}}<\delta$.  Such an
$m$ exists because the first of these two requirements is satisfied by
infinitely many $m$ and the second by all sufficiently large $m$.
Then, adjoining one more point to the domain of $f$ and extending $f$
to take the value $m$ there, we get a condition
$(f\cup(\dom f,m),A,(\eps+\eta)/2)$ that extends $(f,A,\eps)$ and has
the desired inequality.  This completes our preliminary step, and we
assume from now on that
$r\neq\sum_{i<d}\sum_{n\in\dom f}s_ia^i_{f(n)}$.  We introduce the
notation
\[
\zeta=\big|r-\sum_{i<d}\sum_{n\in\dom f}s_ia^i_{f(n)}\big|,
\]
so that we have arranged $\zeta>0$.

As in previous proofs, Corollary~\ref{ls:cor} provides a permutation
$p$ of \bbb N, extending $f$ and satisfying
\[
\sum_n\sq{a^i_{p(n)}:i<d}=\sum_{n\in\dom f}\sq{a^i_{f(n)}:i<d}
\]
and so
\[
\sum_{n\geq\dom f}\sq{a^i_{p(n)}:i<d}=0.
\]
As before, let $\eta<\eps$ be such that $\nm{\sq{a^i_m:i<d}}<\eta/C_d$
for all $m\notin\ran f$, and let $\delta$ be a positive rational
number smaller than both $(\eps-\eta)/2$ and
$\zeta/(2\sum_{i<d}|s_i|)$.  Continuing as in earlier proofs, fix
$n_*$ so large that
\[
\nmb{\sum_{m\in n_*\setminus\dom f}\sq{a^i_{p(m)}:i<d}}<\delta
\]
and $\nm{\sq{a^i_m:i<d}}<\delta/C_d$ for each $m\geq n_*$. By
Corollary~\ref{poly:cor}, there is an injection $g:n_*\to\bbb N$
extending $f$, having range $p[n_*]$, and satisfying
\[
\nmb{\sq{a^i_{\dom f}:i<d}+\sum_{\dom f<k<m}\sq{a^i_{g(k)}:i<d}}
\leq (\eta/C_d)C_d+\delta<\eps-\delta
\]
for all $m\leq n_*$.  As before, this ensures that $(g,A,\delta)$ is
an extension of $(f,A,\eps)$ in $\bbb P_I$.

Finally, comparing sums over $\dom g$ to sums over $\dom f$, we have
that
\[
\nmb{\sum_{m\in\dom g}\sq{a^i_m:i<d}-\sum_{m\in\dom f}\sq{a^i_m:i<d}}
<\delta
\]
and so
\[
\big|\sum_{m\in\dom g}\sum_{i<d}s_ia^i_m-
\sum_{m\in\dom f}\sum_{i<d}s_ia^i_m\big|
<\delta\sum_i|s_i|<\frac\zeta2.
\]
Combining this with the definition of $\zeta$, we find that
\[
\big|r-\sum_{m\in\dom g}\sum_{i<d}s_ia^i_m\big|>
\frac\zeta2>\delta\sum_{i<d}|s_i|,
\]
as required.
\end{pf}

Putting the lemmas together, we obtain the following theorem
describing what forcing by $\bbb P_I$ accomplishes.

\begin{thm}     \label{force1:thm}
Let $I$ be a maximal independent family of conditionally convergent
real  series, and let $G\subseteq\bbb P_I$ be a $V$-generic filter.
Let $\pi=\bigcup\{f:(f,A,\eps)\in G\}$ and let $b$ be any conditionally
convergent series in $V$.  Then $\sum_nb_{\pi(n)}$ converges to a sum
not in $V$.
\end{thm}

\begin{pf}
We have already seen in Corollary~\ref{perm-conv:cor} that $\pi$ is a
permutation of \bbb N and that the rearranged series
$\sum_na_{\pi(n)}$ converges for each $a\in I$.

Because $I$ is a maximal independent set, any conditionally convergent
series $b$ is the sum of an absolutely convergent series $c$ and a linear
combination $\sum_{i<d}s_ia^i$ of some elements $a^i$ of $I$ with
nonzero coefficients.  It follows immediately that $\sum_nb_{\pi(n)}$
converges.  Furthermore, the absolutely convergent $c$ has the same
sum after rearrangement as before; in particular, the rearranged sum
is in $V$.  So to complete the proof of the theorem, it suffices to
show that  $\sum_{i<d}\sum_ns_ia^i_{\pi(n)}$ is not in $V$.  To prove this, we
fix an arbitrary real $r\in V$ and show that
$\sum_{i<d}\sum_ns_ia^i_{\pi(n)}\neq r$.

By Lemma~\ref{dense3} and genericity, $G$ contains a condition
$(g,A,\delta)$ satisfying the conclusion of that lemma.  So
\[
\eta=\big|r-\sum_{i<d}\sum_{n\in\dom g}s_ia^i_{g(n)}\big|
>\delta\sum_{i<d}|s_i|.
\]
Now consider any extension $(h,B,\gamma)\leq(g,A,\delta)$ in
  $G$. We have
\[
\nmb{\sum_{n\in\dom h}\sq{a^i_{h(n)}:i<d}
-\sum_{n\in\dom g}\sq{a^i_{g(n)}:i<d}}<\delta,
\]
so
\[
\big|\sum_{i<d}\sum_{n\in\dom h}s_ia^i_{h(n)}-
\sum_{i<d}\sum_{n\in\dom g}s_ia^i_{g(n)}\big|<\delta\sum_{i<d}|s_i|,
\]
and therefore
\[
\big|r-\sum_{i<d}\sum_{n\in\dom h}s_ia^i_{h(n)}\big|>\eta-
\delta\sum_{i<d}|s_i|>0.
\]

Because the generic filter $G$ is directed, we know that, among the
partial sums of the rearranged series
$\sum_n\sum_{i<d}s_ia^i_{\pi(n)}$, cofinally many are of the form
$\sum_{i<d}\sum_{n\in\dom h}s_ia^i_{h(n)}$ for some $(h,B,\gamma)$ as
above. These partial sums therefore differ from $r$ by more than the
positive constant $\eta-\delta\sum_{i<d}|s_i|$.  Note that this
constant is independent of $(h,B,\gamma)$. We therefore conclude that
the infinite sum $\sum_n\sum_{i<d}s_ia^i_{h(n)}$ differs from
$r$ by at least $\eta-\delta\sum_{i<d}|s_i|$ and is therefore
certainly not equal to $r$.
\end{pf}

In order to iterate forcings of the form $\bbb P_I$, we use the chain
condition provided by the following lemma.

\begin{la}      \label{ccc1}
$\bbb P_I$ satisfies the countable chain condition.  
\end{la}

\begin{pf}
Let $\scr A$ be an uncountable set of conditions in $\bbb P_I$. By Lemma~\ref{dense1}, we may for each condition in $\scr A$ find a stronger condition $(f,A,\eps)$ such that $\nm{\sq{a^i_m:i<d}} < \eps/2C_{2d}$ for all $m \notin \mathrm{Range}(f)$, where $d = |A|$ and $\sq{a^i:i<d}$ is an enumeration of $A$. Let $\scr A'$ denote some set obtained from $\scr A$ by replacing each condition with a stronger condition in this way. If $\scr A'$ is countable then we are done, so let us suppose $\scr A'$ is uncountable.

It is straightforward to verify that, for two conditions $(f,A,\eps)$
and $(g,B,\delta)$ in $\bbb P_I$ to be compatible, it is sufficient to have
\begin{ls}
  \item $f=g$,
\item $\eps=\delta$,
\item $|A|=|B|$, and
\item if $A$ and $B$ are enumerated as \sq{a^i:i<d} and \sq{b^i:i<d},
  then \nm{\sq{a^i_m:i<d}} and \nm{\sq{b^i_m:i<d}} are each
  $<\eps/2C_{2d}$ for all $m\notin\ran f$.  
\end{ls}
Since there are only
  countably many possibilities for $f$, for $\eps$ (recall that $\eps$
  has to be rational), and for $|A|$, it follows that some two conditions in $\scr A'$ must be compatible. But this shows that some two conditions in $\scr A$ must be compatible.
\end{pf}

Combining this lemma with Theorem~\ref{force1:thm}, we obtain the
following result for a finite-support iteration.

\begin{thm}
Suppose $\ger c>\aleph_1$.  Let $\bbb P$ be a finite-support iteration
of length $\omega_1$ where each stage of the forcing is $\bbb P_I$ for
some maximal independent set of conditionally convergent series in the
extension produced by the previous stages of the iteration.  Then in
the extension produced by forcing with \bbb P, we have
$\rrf=\aleph_1<\ger c$.
\end{thm}

\section{Forcing Infinite Limits}       \label{force2}

In this section, we describe a notion of forcing \bbb P  producing
a permutation $\pi$ that rearranges all conditionally convergent
series in the ground model so that they diverge to $+\infty$ or to
$-\infty$.  Afterward, we iterate this forcing and the one from the
previous section to show that all our rearrangement numbers can
consistently be strictly smaller than \ger c.

In fact, the argument here applies not only to conditionally
convergent series but to a broader class of series defined as follows.

\begin{df}      \label{pcc:df}
A series of real numbers is \emph{potentially conditionally
  convergent,} abbreviated \emph{pcc,} if some rearrangement of it is
conditionally convergent.
\end{df}

It is easy to see that a series is pcc if and only if its terms
converge to zero and the two sub-series consisting of its positive
terms and its negative terms both diverge.

\begin{rmk}
  Readers who are interested only in conditionally convergent series,
  not in pcc ones, can safely interpret ``pcc'' in the rest of this
  section as meaning conditionally convergent.  Another safe
  simplification of most of the the following material (all but
  Corollary~\ref{sigmacent:cor})  is that readers
  uncomfortable with the version $\text{MA}(\sigma\text{-centered})$
  of Martin's Axiom used below can pretend that we refer to the
  ordinary, stronger version MA.
\end{rmk}

\begin{df}
Let $\bar a=\sq{a_n:n\in\bbb N}$ be a sequence of real numbers.
Define
\begin{ls}
\item $P(\bar a) =\{n\in\bbb N:a_n>0\}$,
\item $N(\bar a) =\{n\in\bbb N:a_n<0\}$,
\item $\scr I(\bar a)=\{A\subseteq\bbb N:\sum_{n\in A}|a_n|\text{
    converges}\}$,
\item $\scr I^+(\bar a)=\{A\subseteq\bbb N:\sum_{n\in A}|a_n|\text{
    diverges}\}$,  and
\item $\scr I^*(\bar a)=\{A\subseteq\bbb N:\bbb N\setminus A\in\scr
  I(\bar a)\}$.
\end{ls}
\end{df}

$\scr I(\bar a)$ is known in the literature as the \emph{summability
  ideal} for $\bar a$ (or, more precisely, for the sequence of absolute
values \sq{|a_n|:n\in\bbb N}).  The terminology ``ideal'' is justified
because $\scr I(\bar a)$ is clearly closed under subsets and under
finite unions.  Its complement $\scr I^+(\bar a)$ is the associated
co-ideal and $\scr I^*(\bar a)$ is the associated filter.

\begin{notat}
In preparation for defining the desired forcing \bbb P, we fix an
enumeration, of length \ger c, of all the pcc series in the ground
model.  We regard each series $\sum_na_n$ as the sequence
\sq{a_n:n\in\bbb N} of its terms, so we are dealing with a \ger
c-enumeration \sq{\bar a^\beta:\beta<\ger c} of infinite sequences.

We shall also use the standard notation $\subseteq^*$ for
almost-inclusion; that is, $X\subseteq^* Y$ means that $X\setminus Y$
is finite.
\end{notat}

Next, we need a technical lemma.

\begin{rmk}     \label{matrix:rmk}
This remark is intended to clarify the intentions behind
  Lemma~\ref{matrix:la} below.  Of course, in principle, the lemma
  can stand on its own; only the lemma itself, not the intentions,
  will be strictly needed in what follows.

  The lemma is intended to address the same issue already mentioned in
  Remark~\ref{indep:rmk}, namely that correlations between various
  series may constrain our options for dealing with them.  In the
  present situation, it turns out that all the pcc series (in the
  ground model) can be organized into equivalence classes such that
  decisions about one series (for example, whether its
  $\pi$-rearrangement should diverge to $+\infty$ rather than
  $-\infty$) affect the other series in its equivalence class, but do
  not affect series in other equivalence classes.

  The construction of the equivalence classes is complicated by
  the following considerations. If a series $\sum_na_n$ is to be
  rearranged by $\pi$ to diverge to, say, $+\infty$, then $\pi$ must
  move some set, say $X$, of numbers from $P(\bar a)$ to relatively
  earlier positions. We shall want to do this without disturbing
  series $\sum_nb_n$ from other equivalence classes. So it is
  desirable that the moved numbers from $P(\bar a)$ constitute a set
  $X$ in $\scr I^+(\bar a)$ (so that we can get large partial sums
  this way) but in $\scr I(\bar b)$ (so that $\sum_nb_n$ is not
  seriously disturbed).  So we need that such an $X$ exists. But more
  is needed, because we may have already chosen some set $X'$ of
  numbers to be moved for the sake of some other series $\sum_nc_n$,
  and so we shall need an appropriate $X$ disjoint (or at least almost
  disjoint) from $X'$.  If no such $X$ is available, then we cannot
  handle $\bar b$ independently from $\bar a$, so they will have to go
  into the same equivalence class.

Thus, the choice of the appropriate sets $X$ depends on the equivalence
relation (because elements in the same equivalence class should use
the same $X$'s) but also influences the equivalence relation.  As a
result, the construction of the equivalence classes and the choice of
the $X$'s need to be done in a mutual recursion.  That is what
Lemma~\ref{matrix:la} and its proof are about.

In terms of our fixed enumeration \sq{\bar a^\beta:\beta<\ger c} of
all the pcc series, the equivalence relation described above can be
viewed as an equivalence relation on the set \ger c of indices.  For
each equivalence class, we use its first element (smallest ordinal
number) as a standard representative.  In the notation of the lemma,
$A$ will be the set of these representatives, and $\zeta$ will be the
function sending each ordinal $\beta<\ger c$ to the representative of
its equivalence class.  The $X$'s in the preceding discussion will be
$X$'s in the lemma also, but there is an additional complication as
each equivalence class gets not a single $X$ but an almost decreasing
(modulo finite sets) \ger c-sequence of $X$'s.
\end{rmk}

\begin{la}      \label{matrix:la}
Assume \masc.  There exist a set $A\subseteq\ger c$, a function
$\zeta:\ger c\to A$, and a matrix of sets \sq{X^\beta_\alpha:
  \alpha\in A\text{ and }\alpha\leq\beta<\ger c} with the following
properties for all $\beta<\ger c$:
\begin{lsnum}
\item $\zeta(\beta)\leq\beta$ with equality if and only if $\beta\in A$.
\item If $\alpha\in A$ and $\alpha\leq\beta\leq\beta'$, then
  $X^{\beta'}_\alpha\subseteq^* X^\beta_\alpha$.
\item The sets $X^\beta_\alpha$ for $\alpha\in A\cap(\beta+1)$ are
  almost disjoint, i.e., the intersection of any two distinct ones is
  finite.
\item $X^\beta_{\zeta(\beta)}$ is a subset of $P(\bar a^\beta)$ or of $N(\bar
  a^\beta)$.
\item If $\beta\leq\beta'$ then $X^{\beta'}_{\zeta(\beta)}\in\scr
  I^+(\bar a^\beta)$.
\item If $\alpha\in A$ and $\alpha<\zeta(\beta)$ then
  $X^\beta_\alpha\in\scr I(\bar a^\beta)$.
\item All subsets of $X^\beta_{\zeta(\beta)}$ that belong to $\scr
  I^+(\bar a^{\zeta(\beta)})$ also belong to $\scr I^+(\bar
  a^\beta)$.
\end{lsnum}
\end{la}

\begin{rmk}     \label{matrix2:rmk}
Continuing from Remark~\ref{matrix:rmk}, we comment on the ideas
behind the clauses in this lemma.  We regard two ordinals
$\beta,\beta'<\ger c$ as equivalent if
$\zeta(\beta)=\zeta(\beta')$. By clause~(1), $\zeta(\beta)$ is the
first element of the equivalence class of $\beta$, and $A$ is the set
of all these first elements, for all the equivalence classes.
Clause~(7) describes the effect of equivalence of ordinals on the
associated series. Roughly speaking, it correlates divergence of
subseries of $\bar a^\beta$ with divergence of the corresponding
subseries of $\bar a^{\zeta(\beta)}$.  Clauses~(5) and (6) act in the reverse
direction for some (not all)  inequivalent ordinals. Specifically, if
$\alpha<\zeta(\beta)$ are the first elements of two equivalence
classes then $X^\beta_\alpha$ is in $\scr I^+(\bar a^\alpha)$ (by (5)
with $\alpha$ and $\beta$ in place of $\beta$ and $\beta'$) but in
$\scr I(\bar a^\beta)$ (by (6)).

Clause~(4) implies that the matrix of $X$'s decides a direction,
positive or negative, for each pcc series $\bar a^\beta$.  This
decision will later determine whether the generic rearrangement
of $\bar a^\beta$ will diverge to $+\infty$ or to $-\infty$.

Clauses~(2) and (3) describe the general structure of the $X$ matrix.
If we regard the subscripts as the horizontal coordinate and the
superscripts as vertical, then (2) says that the columns are almost
decreasing, and (3) says that the rows are almost disjoint.
\end{rmk}

We now turn to the proof of the lemma.

\begin{pf}
We proceed by recursion on ordinals $\beta<\ger c$.  At stage $\beta$,
we shall define $\zeta(\beta)$, $A\cap(\beta+1)$, and the $\beta\th$
row \sq{X^\beta_\alpha:\alpha\in A\cap(\beta+1)} of the $X$ matrix.

For $\beta=0$, we set $\zeta(0)=0$ (as required by clause~(1) of the
lemma) and we put 0 into $A$ (as required by $\zeta:\ger c\to A$). For
$X^0_0$, we must take a set that is in $\scr I^+(\bar a^0)$ (as
required by (5)) and that is a subset of $P(\bar a^0)$ or of
$N(\bar a^0)$ as required by (4).  Such sets exist, i.e., $\bar a^0$
has divergent subseries consisting of only positive terms or only
negative terms, because $\bar a^0$ is pcc.  (In fact, we can choose
$P$ or $N$ here as we wish; both sorts of sets exist.)

Next, we consider the case of successor ordinals.  Suppose stage
$\beta$ has been completed, so, in particular, we have the almost
disjoint family $\{X^\beta_\alpha:\alpha\in A\cap(\beta+1)\}$.  To
produce the required items for $\beta+1$, we proceed by a subsidiary
recursion on $\alpha\in A\cap(\beta+1)$ as follows.

At step $\alpha$ of this recursion, we consider two cases, according
to whether or not there exists a subset $Y$ of $X^\beta_\alpha$ that is
in $\scr I^+(\bar a^\alpha)\cap\scr I(\bar a^{\beta+1})$.

If such a $Y$ exists, then we choose one and declare it to be
$X^{\beta+1}_\alpha$.  Then we proceed to the next value of $\alpha$.

If no such $Y$ exists, then we stop the subsidiary recursion on
$\alpha$, we define $\zeta(\beta+1)=\alpha$, and we declare
$\beta+1\notin A$ (as required by clause~(1)), so
$A\cap(\beta+2)=A\cap(\beta+1)$. We define $X^{\beta+1}_\alpha$ to be
some subset of $X^\beta_\alpha$ that is included in either $P(\bar
a^{\beta+1})$ or $N(\bar a^{\beta+1})$.  To see that such a set
exists, notice first that $X^\beta_\alpha$ is in $\scr I^+(\bar
a^\alpha)$ because the earlier stage $\beta$ of our main recursion
satisfied clause~(5) (and $\alpha\in A$).  Next, use the case
hypothesis to infer that $X^\beta_\alpha\in\scr I^+(a^{\beta+1})$,
which means that the series $\sum_{n\in X^\beta_\alpha}|a^{\beta+1}_n|$
diverges. Finally infer that, in this divergent series, either the
positive terms or the negative terms form a divergent series, and the
index set of such a series can serve as the desired
$X^{\beta+1}_\alpha$.  Finally, we set
$X^{\beta+1}_\gamma=X^\beta_\gamma$ for all $\gamma\in A$ in the range
$\alpha<\gamma\leq\beta$.  These choices satisfy all the clauses of
the lemma.

If the subsidiary recursion is not stopped at a stage where no $Y$ is
available, i.e., if this recursion continues through all ordinals in
$A\cap(\beta+1)$, then we have defined $X^{\beta+1}_\alpha$ for all
$\alpha\in A\cap(\beta+1)$; we have not defined $\zeta(\beta+1)$ yet,
nor have we added any element to $A$.  We now put $\beta+1$ into $A$
and, as required by clause~(1), we set $\zeta(\beta+1)=\beta+1$.  We
must still choose a set to serve as $X^{\beta+1}_{\beta+1}$.  This set
must be
\begin{ls}
\item in $\scr I^+(\bar a^{\beta+1})$ (by clause~(5)),
\item  almost disjoint from all the sets $X^{\beta+1}_\alpha$ for
$\alpha\in A\cap(\beta+1)$ (by clause~(3)), and
\item a subset of $P(\bar a^{\beta+1})$ or of $N(\bar a^{\beta+1})$
  (by clause~(4)).
\end{ls}
If we can find such a set then, by using it as
$X^{\beta+1}_{\beta+1}$, we shall satisfy all the clauses for this
stage $\beta+1$.  Furthermore, any set satisfying the first two of
these three requirements can be pruned to satisfy the third.  This is
because, as noted before, if a series diverges then either the
subseries of positive terms or the subseries of negative terms (or
both) will also diverge.

So to complete the successor stage of our induction on $\beta$, we
must prove the existence of a set in the co-ideal $\scr I^+(\bar
a^{\beta+1})$ that is almost disjoint from all the sets $X^{\beta+1}_\alpha$ for
$\alpha\in A\cap(\beta+1)$.  It is here that we must invoke \masc.

Specifically, we apply \masc\ to Mathias forcing guided by the filter
$\scr I^*(\bar a^{\beta+1})$.  Forcing conditions are pairs $(s,C)$
where $s$ is a finite subset of \bbb N and
$C\in\scr I^*(\bar a^{\beta+1})$ with $\min(C)>\max(s)$. Another
condition $(s',C')$ is an extension of $(s,C)$ if $s\subseteq s'$,
$C\supseteq C'$, and $s'\setminus s\subseteq C$. This forcing is
$\sigma$-centered (and thus satisfies the countable chain condition)
because any finitely many conditions with the same first component are
compatible; just intersect their second components. So we can apply
\masc\ with the following fewer than \ger c dense sets.

First, for each of the sets $X^{\beta+1}_\alpha$ that we want our
$X^{\beta+1}_{\beta+1}$ to be almost disjoint from, we have the dense
set
\[
D_\alpha=\{(s,C):X^{\beta+1}_\alpha\cap C=\emp\}.
\]
There are fewer than \ger c of these sets, as they are indexed by
ordinals $\alpha\in A\cap(\beta+1)$, and each of them is dense because
the sets $X^{\beta+1}_\alpha$ were chosen, in our subsidiary
recursion, to be in the ideal $\scr I(\bar a^{\beta+1})$.

Second, for each natural number $k$, we have the dense set
\[
D'_k=\{(s,C):|\sum_{n\in s}a^{\beta+1}_n|>k\}.
\]
This is dense because the second components of our conditions are sets
in $\scr I^*(\bar a^{\beta+1})$ and the summation of $\bar
a^{\beta+1}$ over any such set is pcc.

By \masc, there is a filter $G$ of conditions meeting all these dense
sets.  Let $X^{\beta+1}_{\beta+1}=\bigcup\{s:(s,C)\in G\}$. The sum of
$\bar a^{\beta+1}$ over $X^{\beta+1}_{\beta+1}$ diverges because $G$
meets every $D'_k$.  And the fact that $X^{\beta+1}_{\beta+1}$ is
almost disjoint from each previous $X^{\beta+1}_\alpha$ follows by a
routine compatibility argument from the fact that $G$ meets every
$D_\alpha$.

This completes the recursion for successor steps $\beta+1$. We turn to
the limit case.

Let $\beta$ be a limit ordinal, and suppose the construction has been
carried out, in accordance with the requirements of the lemma, for all
$\gamma<\beta$.  For each $\alpha\in A\cap\beta$, we shall first
produce a set $Y_\alpha\in\scr I^+(\bar a^\alpha)$ such that
$Y_\alpha\subseteq^* X^\gamma_\alpha$ for all $\gamma<\beta$.  Once
this is done, we can proceed exactly as in the successor case, using
$Y_\alpha$ in place of $X^{\beta}_\alpha$ and defining sets called
$X^\beta_\alpha$ rather than $X^{\beta+1}_\alpha$.

To produce the desired $Y_\alpha$, we consider any fixed
$\alpha\in A\cap\beta$ and apply \masc\ to Mathias forcing guided by
the filter generated by $\scr I^*(\bar a^\alpha)$ and the sets
$X^\gamma_\alpha$ for $\alpha\leq\gamma<\beta$.  This is a proper
filter because the sets $X^\gamma_\alpha$ that we are adjoining to
$\scr I^*(\bar a^\alpha)$ form an almost decreasing sequence (by
clause~(2) of the lemma for stages $\gamma<\beta$) of sets in
$\scr I^+(\bar a^\alpha)$ (by clause~(5)). The relevant dense sets are
\[
D_\gamma=\{(s,C):C\subseteq X^\gamma_\alpha\}
\]
(dense because $X^\gamma_\alpha$ is in the guiding filter) and
\[
D'_k=\{(s,C):|\sum_{n\in s}a^{\beta}_n|>k\}
\]
as in the earlier use of \masc.  A generic filter $G$ meeting all
these dense sets produces the desired $Y_\alpha=\bigcup\{s:(s,C)\in
G\}$.  As before $\sum_{n\in Y_\alpha}a^\beta_n$ diverges because $G$
meets the dense sets $D'_k$, and $Y_\alpha\subseteq^* X^\gamma_\alpha$
because $G$ meets $D_\gamma$.

This completes the proof of the existence of the required
$Y_\alpha$'s, and thus completes the recursion on $\beta$ that
produces the sets and function required by the lemma.
\end{pf}

Fix $A$, $\zeta$, and \sq{X^\beta_\alpha} as in the lemma.  Call an
ordinal $\beta<\ger c$ a P-ordinal or an N-ordinal according to whether
the P or N alternative holds in clause~(4) of the lemma.
(Intuitively, the P-ordinals are those for which the generic
rearrangement of $\sum_na^\beta_n$ will diverge to $+\infty$, and the
N-ordinals are those for which this rearrangement will diverge to
$-\infty$.)  We write $R(\beta)$ for $P(\bar a^\beta)$ when $\beta$ is
a P-ordinal and for $N(\bar a^\beta)$ when $\beta$ is an N-ordinal.

Let \bbb P be the following forcing.  A condition is a triple
$(s,F,k)$ such that
\begin{ls}
  \item $s$ is an injective function from some $n\in\bbb N$ into \bbb
    N,
\item $F$ is a finite subset of \ger c,
\item $k\in\bbb N$, and
\item For all P-ordinals (resp.\ N-ordinals) $\beta\in F$, the sum
  $\sum_{i\in\dom s}a^\beta_{s(i)}$ is positive (resp.\ negative) and its
  absolute value is $>k$.
\end{ls}
A condition $(s',F',k')$ extends $(s,F,k)$ if
\begin{ls}
  \item $s'\supseteq s$,
\item $F'\supseteq F$,
\item $k'\geq k$, and
\item for all $j\in\dom{s'}\setminus\dom s$ and all P-ordinals (resp.\
  N-ordinals) $\beta\in F$, $\sum_{i<j}a^\beta_{s'(i)}$ is positive (resp.\
  negative) and its absolute value is $>k$.
\end{ls}

Intuitively, the intended ``meaning'' of a condition $(s,F,k)$ is that
the generic permutation $\pi$ will be an extension of $s$ and that the
finitely many series $\sum_na^\beta_{\pi(n)}$ for $\beta\in F$ are
well on their way to diverging in the intended direction.  Here
``well'' is measured by $k$, and ``well on their way'' means that the
partial sum provided by $s$, namely $\sum_na^\beta_{s(n)}$, exceeds $k$
in the intended (positive or negative) direction and will continue to
do so as more terms are included in the partial sum.  That
$\sum_na^\beta_{s(n)}$ is large enough in the appropriate direction is
the content of the last clause in the definition of conditions; that
longer partial sums will also behave in this way is the content of the
last clause in the definition of extensions.

\begin{la}      \label{dense4:la}
  Assume \masc.
\begin{lsnum}
\item \bbb P satisfies the countable chain condition; in fact, it is
  $\sigma$-centered.
\item For every $l\in\bbb N$, every condition $(s,F,k)$ has an
  extension $(s',F',k')$ with $k'\geq l$.
\item For every $m\in\bbb N$, every condition $(s,F,k)$ has an
  extension $(s',F',k')$ with $m\in\ran{s'}$.
\item For every $\gamma\in\ger c$, every condition $(s,F,k)$ has an
  extension $(s',F',k')$ with $\gamma\in F'$.
  \end{lsnum}
\end{la}

\begin{pf}
Part~(1) is easy. Any finitely many conditions $(s,F,k)$ with the same
$s$ and $k$ are compatible; just take the union of the $F$'s.

For part~(2), it suffices to treat the case $l=k+1$, since repeated
extensions of this sort yield arbitrarily large $l$'s.  Let $(s,F,k)$ be
given; the desired extension will be of the form $(s',F,k+1)$ (with
the same $F$). Our task is to produce an $s'\supseteq s$ such that
this $(s',F,k+1)$ is an extension of $(s,F,k)$, and that comes down to
satisfying the last clause in the definition of condition and the last
clause in the definition of extension.

Before proceeding with the detailed proof, we describe the idea behind
it; this paragraph can be omitted by readers who just want the
detailed proof. We shall extend $s'$ in several steps, where each step
serves to make the partial sum $\sum_{n\in\dom{s'}}a^\beta_{s'(n)}$
for some $\beta\in F$ appropriately large; these sums, taken only up
to $\dom s$, were already bigger than $k$ in absolute value; the
extension $s'$ must make them bigger than $k+1$. The difficulty is
that, when we extend $s'$ to make one of these sums, say the one for
$\beta$, large, there is a danger of making other sums, for other
$\beta'\in F$, too small, and we cannot afford to do this. Not only
must the final sums, over all of $\dom{s'}$ be at least $k+1$, but
they cannot drop below $k$ at any point between \dom s and \dom{s'}
(because of the last clause in the definition of extension).  This
difficulty will be overcome by means of two observations.  First, if
$\zeta(\beta)<\zeta(\beta')$, then clauses~(5) and (6) of
Lemma~\ref{matrix:la} provide a set
$X^{\max\{\beta,\beta'\}}_{\zeta(\beta)}$ on which the series
$\sum a^\beta_n$ diverges while $\sum a^{\beta'}_n$ converges. This
means that we can append finitely many elements of that set to the
range of $s$ in such a way as to get a big (in absolute value) partial
sum for the $\beta$ series while making very little change in the
$\beta'$ series.  The partial sum for $\beta'$ may get smaller, but
not too small.  In other words, when we want to make the partial sum
of $a^\beta$ large, we need not worry about causing trouble for
$\beta'$ with larger $\zeta$ values.  What about $\beta'$ with smaller
$\zeta$ values? They might get seriously damaged by what we do for
$\beta$, but one could recover from that damage by making the partial
sum of $\bar a^{\beta'}$ very large, not just bigger than $k+1$ (our
original goal) but so much bigger that the damage from $\beta$ is
cancelled.  We thus overcompensate at $\beta'$ for the damage done by
$\beta$. This does not quite solve the problem, because we need to
control not only the final partial sums over all of \dom{s'} but also
the intermediate partial sums for $j$ as in the last clause of the
definition of extension.  This, fortunately, can be easily handled: We
do the overcompensation before the damage.  By putting terms into $s'$
in the right order, we can first make the partial sums for $\beta'$
very large (overcompensation), and have the damage come later so that
the partial sums stay large all the time.

Here are the formal details implementing the ideas in the preceding
paragraph.  Let $\delta$ be the largest element of $F$. (If $F$ is
empty, the construction is trivial.) Enumerate $\zeta[F]$ in
increasing order as $\zeta_0<\zeta_1<\dots<\zeta_{m-1}$.  Since $(s,F,k)$
is a condition, and since the last requirement in the definition of
conditions demanded a strict inequality, fix an $\eps>0$ such that,
for all $\beta\in F$, $\sum_{n\in\dom s}|a^\beta_n|>k+\eps$.

By backward recursion on $j<m$, find finite sets $Z_j\subseteq\bbb N$
such that
\begin{ls}
\item the sets $Z_j$ are disjoint from each other and from \ran s,
\item $Z_j\subseteq X^\delta_{\zeta_j}\cap\bigcap\{R(\beta): \beta\in
  F\text{ and }\zeta(\beta)=\zeta_j\}$,
\item for all $\beta\in F$ with $\zeta(\beta)=\zeta_j$,
\[
\big|\sum_{i\in Z_j}a^\beta_i\big|>1+\sum_{j'>j}\sum_{i\in Z_{j'}}|a^\beta_i|,
\]
\item for all $\beta\in F$ with $\zeta(\beta)>\zeta_j$,
\[
\sum_{i\in Z_j}|a^\beta_i|<\frac\eps m.
\]
\end{ls}
To see that such sets $Z_j$ can be chosen, consider a particular $j$
and suppose appropriate $Z_{j'}$ have already been chosen for all $j'$
in the range $j<j'<m$. According to the second requirement, we seek
$Z_j$ as a subset of $X^\delta_{\zeta_j}$.  The rest of the second
requirement prohibits only finitely many elements of this set from
being in $Z_j$, thanks to clauses~(2) and (4) of Lemma~\ref{matrix:la}
(remember that $\delta\geq\beta$ for all $\beta\in F$).  For the first
requirement, we exclude finitely many more elements from potentially
entering $Z_j$, namely the elements of \ran s and the elements of the
previously chosen $Z_{j'}$ for $j'>j$.  The fourth requirement also
excludes only finitely many elements of $X^\delta_{\zeta_j}$, because
the relevant series converge absolutely when restricted to
$X^\delta_{\zeta_j}$, thanks to clause~(6) of Lemma~\ref{matrix:la}.
So we have a cofinite subset of $X^\delta_{\zeta_j}$ in which to find
$Z_j$ satisfying the third requirement.  And this task is easy because
the series $\sum_ia^\beta_i$ restricted to this set diverges by
clause~(5) of Lemma~\ref{matrix:la} (and has all its terms of the same
sign by the second requirement).

Now that we have the $Z_j$'s, we use them to define $s'$ as follows.
It is the extension of $s$ obtained by appending the elements of all
the $Z_j$'s in order of increasing $j$.  More formally, the domain of
$s'$ is $\dom s+\sum_{j<m}|Z_j|$ and
\[
s'(i)=
\begin{cases}
  s(i)&\text{if }i\in\dom s,\\
  t\th\text{ element of }Z_j& \text{if }i=\dom s+\sum_{j'<j}|Z_{j'}|+t
  \text{ and }t<|Z_j|.
\end{cases}
\]
To check that $(s',F,k+1)$ is a condition extending
$(s,F,k)$, as required for part~(2) of the lemma, we need only check
the last clause in the definition of condition and the last clause in
the definition  of extension.  So we need to consider sums of the form
$\sum_{n\in\dom t}a^\beta_{t(n)}$ where $\beta\in F$ and $t$ is an
initial segment of $s'$ that strictly includes $s$. We need these sums to be
$>k$ for all such $t$ and $> k+1$ when $t$ is all of $s'$.

For brevity, we shall say ``block $j$ of $s'$'' for the segment of
$s'$ that corresponds to the $j$ part of the definition above, that
is, $s'$ restricted to the interval
\[
[\dom s+\sum_{j'<j}|Z_j|,\dom s+\sum_{j'\leq j}|Z_j|);
\]
this is the part of $s'$ whose range is $Z_j$.

Fix some $\beta\in F$ and let $q$ be the index such that
$\zeta(\beta)=\zeta_q$. For convenience, assume $\beta$ is a
P-ordinal; the proof for N-ordinals is the same up to signs.

Consider the sums $S(x)=\sum_{n<x}a^\beta_{s'(n)}$ beginning with $x=\dom
s$ and gradually increasing $x$ up to \dom{s'}.  We begin with $S(\dom
s)>k+\eps$ by our choice of $\eps$ (ultimately coming from the fact
that $(s,F,k)$ is a condition). As $x$ increases through blocks of
$s'$ that strictly precede block $q$, $S(x)$ changes by less than
$\eps/m$ in any single block, by the last clause in our choice of the
$Z_j$'s.  There are only $m$ blocks altogether, so the total variation
in $S(x)$ before block $q$ is less than $\eps$. Since $S(\dom s)$ was
$>k+\eps$, this variation cannot bring $S(x)$ down to $k$ or lower.

Continuing to increase $x$, into block $q$, we find only positive
summands, by the second clause in our choice of $Z_j$'s, so $S(x)$
increases and, in particular, remains $>k$, while $x$ is in block
$q$.  By the end of block $q$, $S(x)$ has become quite large, thanks
to clause 3 in our choice of $Z_j$'s.  Specifically, using clause~3
and the fact that $S(x)$ was $>k$ at the beginning of block $q$, we
see that, at the end of block $q$, $S(x)$ has grown to more than
\[
k+1+\sum_{j'>q}\sum_{i\in Z_{j'}}|a^\beta_i|.
\]
In this formula, the sum over $j'$ and $i$ majorizes the absolute
value of any further change in $S(x)$ beyond block $q$.  So, from the
end of block $q$ on, $S(x)$ will always be more than $k+1$.  This
ensures that $(s',F,k+1)$ is a condition and also that it is an
extension of $(s,F,k)$, as required.  This completes the proof of
part~(2) of the lemma.

Part~(3) is a consequence of part~(2) as follows. Let $(s,F,k)$ be any
condition, and let $m\in\bbb N$. The desired result is trivial if
$m\in\ran s$, so assume $m\notin\ran s$.  (Intuition: We cannot simply
append $m$ to the range of $s$, as that might ruin the requirements,
in the definition of condition and extension, that certain partial
sums must remain big in absolute value.  So we first use part~(2) to
make the relevant partial sums big enough that appending $m$ won't
hurt.) Apply part~(2) of the lemma with some $l>k+|a^\beta_m|$ for all
$\beta\in F$.  We get a condition $(s',F,k')$ with $k'\geq l$.  (Recall
that the proof of part~(2) did not require changing $F$.)  Now obtain
$s''$ by appending $m$ to the range of $s'$.  With $F$ and $k$
unchanged (note: $k$, not $k'$), we get that $(s'',F,k)$ is an
extension of $(s,F,k)$ as required.

Finally, to prove part~(4) of the lemma, let any condition $(s,F,k)$
and any ordinal $\gamma<\ger c$ be given; we assume $\gamma\notin F$
because otherwise the conclusion is trivial.  We also assume that
$\gamma$ is a P-ordinal; the proof for an N-ordinal is the same except
for some minus signs.  It suffices to produce an extension of
$(s,F,k)$ of the form $(s',F\cup\{\gamma\},k)$.  (Intuition: We cannot
in general take $s'=s$, because $\sum_{i\in\dom s}a^\gamma_{s(i)}$
might be smaller than $k$, and then $(s,F\cup\{\gamma\},k)$ would fail
to satsify the last requirement in the definition of condition.  We
must extend $s$ to an $s'$ that makes $\sum_{i\in\dom
  s}a^\gamma_{s'(i)}$ large enough.  But we must ensure that we do
not, in this extension, ruin the largeness of the sums for ordinals
$\beta\in F$.  The strategy for doing this is essentially the same as
in the proof of part~(2) above.)

As in the proof of part~(2), let $\delta$ be the largest element of
$F\cup\{\gamma\}$.  Enumerate the ordinals in
$\zeta[F]\cap\zeta(\gamma)$ (not all of $\zeta[F]$) in increasing
order as $\zeta_0<\zeta_1<\dots<\zeta_{m-1}$.  Abbreviate
$\zeta(\gamma)$ as $\zeta_m$, so we have
$\zeta_0<\zeta_1<\dots<\zeta_{m-1}<\zeta_m$.  By backward recursion
on $j\leq m$, find finite sets $Z_j\subseteq\bbb N$ such that
\begin{ls}
\item the sets $Z_j$ are disjoint from each other and from \ran s.
\item $Z_j\subseteq X^\delta_{\zeta_j}\cap\bigcap\{R(\beta):\beta\in
  F\cup\{\gamma\} \text{ and }\zeta(\beta)=\zeta_j\}$,
\item For all $\beta\in F$ with $\zeta(\beta)=\zeta_j$,
\[
\big|\sum_{i\in Z_j}a^\beta_i\big|>\sum_{j'>j}\sum_{i\in
  Z_{j'}}|a^\beta_i|,
\]
\item \[
\big|\sum_{i\in Z_m}a^\gamma_i\big|>k+\big|\sum_{i<\dom
  s}a^\gamma_{s(i)}\big|,
\]
\item for all $\beta\in F\cup\{\gamma\}$ with $\zeta(\beta)>\zeta_j$,
\[
\sum_{i\in Z_j}|a^\beta_i|<\frac\eps m.
\]
\end{ls}
These requirements resemble those in the proof of part~(2).  The
differences are that $j$ ranges up to and including $m$, that $\gamma$
is included along with the elements of $F$ in the second and fifth
clauses, that the right side of the inequality in the third clause
doesn't need an added 1, and that there is a fourth clause
specifically about $\gamma$.  That fourth clause is similar in spirit
to the third, but it takes into account that we have no information
about $\sum_{i\in\dom s}a^\gamma_{s(i)}$; in particular this sum might
be a large negative number.  The right side of the inequality in the
fourth clause is designed to compensate for any such negativity and to
add $k$ beyond that.

The same argument as in the proof of part~(2) yields the existence of
sets $Z_j$ satisfying these requirements.  Once we have these sets, we
can define $s'$ just as we did for part~(2), except of course that now
$j$ ranges up to and including $m$ where we previously had $j<m$.  The
proof that $(s',F\cup\{\gamma\},k)$ is a condition extending $(s,F,k)$
is also just as it was in part~(2).  This completes the proof of
Lemma~\ref{dense4:la}.
\end{pf}

Assume \masc, and suppose $G\subseteq\bbb P$ is a $V$-generic filter.
Let $\pi=\bigcup\{s:(s,F,k)\in G\}$.  Since the first components $s$ of
conditions $(s,F,k)$ are injective functions $n\to\bbb N$ with
$n\in\bbb N$, and since compatible conditions, such as those in $G$,
have compatible first components, $\pi$ is an injective function from
an initial segment of \bbb N into \bbb N.  By part~(3) of
Lemma~\ref{dense4:la} and genericity, $\pi$ is surjective, and so its
domain must be all of \bbb N.  Thus, $\pi$ is a permutation of \bbb
N.  Parts~(2) and (4) of Lemma~\ref{dense4:la}, genericity, and the
last clause in the definition of extensions in \bbb P ensure that
$\sum_na^\beta_{\pi(n)}$ diverges to $+\infty$ for all P-ordinals
$\beta$ and diverges to $-\infty$ for all N-ordinals $\beta$.  So we
have a $\sigma$-centered forcing producing a permutation that
rearranges all pcc series in the ground model to diverge to $+\infty$
or to $-\infty$.

We can achieve the same result without assuming \masc\ in the ground
model.  Simply use, instead of \bbb P, a two-step iteration in which
the first step forces \masc\ and the second step is \bbb P.

Iterating such a forcing for $\omega_1$ steps with finite supports
over a model of $\ger c>\aleph_1$ produces a model in which
$\rri=\aleph_1<\ger c$.  Indeed, the $\aleph_1$ permutations $\pi$
adjoined by the steps of the iteration witness that $\rri=\aleph_1$
because, thanks to the countable chain condition,  every conditionally
convergent series appears
at some intermediate stage of the iteration and is rearranged to
diverge to $+\infty$ or $-\infty$ by the $\pi$ added at the next step
(or at any later step).

We can do even better by combining the forcings from the present
section and Section~\ref{force1}. Start with a model where $\ger
c>\aleph_1$ and perform a finite-support iteration of length
$\omega_1$ in which the steps are alternately the forcing from
Section~\ref{force1} and the forcing from the present section
(including, each time, the forcing of \masc\ that makes \bbb P forcing
possible).  Then each conditionally convergent series in the resulting
model appears at some intermediate stage. The next stage that forces
with Section~\ref{force1}'s forcing (resp.\ the present section's
forcing) produces a permutation making that series converge to a new
finite sum (resp.\ diverge to $+\infty$ or to $-\infty$).  Thus, we
have established the following consistency result:

\begin{thm}     \label{con-small:thm}
  It is consistent with ZFC that $\rrf=\rri=\aleph_1<\ger c$.
\end{thm}

In view of the inequalities between the various rearangement numbers
(see Figure~\ref{fig2}), it follows that all rearrangement numbers can
be $\aleph_1$ while \ger c is larger.  The proof shows that, in
addition, we can make \ger c as large as we wish.

Zapletal \cite{zap} introduced the notion of \emph{tame} cardinal
characteristics of the continuum.  These are characteristics with
definitions of the form ``the smallest cardinality of a set $A$ of
reals such that $\phi(A)\land\psi(A)$,'' where all quantifiers in
$\phi(A)$ range over $A$ or over \bbb N and where $\psi(A)$ has the
form $(\forall x\in\bbb R)(\exists y\in A)\,\theta(x,y)$, where
quantifiers in $\theta(x,y)$ range over \bbb R or \bbb N and where $A$
is not mentioned in $\theta(x,y)$.  By suitable coding (e.g.,
representing a conditionally convergent series by a single real
number), one can show that our rearrangement numbers are tame in
Zapletal's sense; in fact, one doesn't need the $\phi$ component in
the definition of tameness.

Zapletal showed, in \cite[Theorem~0.2]{zap}, that if \ger x is any
tame cardinal characteristic such that $\ger x<\ger c$ holds in some
set-forcing extension, and if there is a proper class of measurable
Woodin cardinals, then $\ger x<\ger c$ holds in the iterated Sacks
model, i.e., the result of a $\ger c^+$-step, countable support
iteration of Sacks forcing.  Thus, we obtain the following corollary
by combining Zapletal's Theorem~0.2 with the tameness of the
rearrangement numbers and the fact that Theorem~\ref{con-small:thm}
was proved by set-forcing.

  \begin{cor}   \label{sacks:cor}
    Assume that there is a proper class of measurable Woodin
    cardinals.  Then the rearrangement numbers are $\aleph_1$ in the
    iterated Sacks model.
  \end{cor}

We note that, if one iterates Sacks forcing beyond $\omega_2$ steps,
it collapses cardinals and, as a result, the iterated Sacks model will
have $\ger c=\aleph_2$.  In contrast, the models we produced to prove
Theorem~\ref{con-small:thm} allowed \ger c to be arbitrarily large.
Also, unlike the proof of Corollary~\ref{sacks:cor}, our proof of
Theorem~\ref{con-small:thm} used no large-cardinal hypotheses.

\begin{rmk}
The proof of Theorem~\ref{con-small:thm} can be easily modified to
obtain other values for the rearrangement numbers.  For example, it is
consistent with ZFC to have
\[
\rr=\rr_{fi}=\rrf=\rri=\aleph_5\quad\text{and}\quad\ger c=\aleph_{17}.
\]
To prove this, begin with a model where $\ger c=\aleph_{17}$ and
iterate the forcings from Sections~\ref{force1} and \ref{force2} for
$\aleph_5$ steps (with finite support as before, and including a
preliminary forcing of \masc\ before each use of the
Section~\ref{force2} forcing).  The same arguments as before show that
the cardinal of the continuum remains $\aleph_{17}$ and that now
$\rrf$ and $\rri$ are at most $\aleph_5$. That not even \rr\ can be
smaller than $\aleph_5$ follows from Theorem~\ref{aboveb:thm}, because
the repeated
forcing of \masc\ makes $\ger b\geq\aleph_5$.
\end{rmk}

In the following corollary, we extract a consistency result from the
fact that the forcing used in the present Section~\ref{force2} is $\sigma$-centered.

\begin{cor}             \label{sigmacent:cor}
It is consistent with ZFC that $\rri<\non N$.
\end{cor}

\begin{pf}
  Consider the forcing with which we obtained $\rri<\ger c$, namely to
  begin with $\ger c>\aleph_1$ in the ground model and to iterate for
  $\omega_1$ steps, with finite supports, the two-step forcing that
  first forces \masc\ and then forces with \bbb P as defined above.
  Let us also suppose that the ground model with which we begin the
  iteration satisfies MA.  We saw above that the final model resulting
  from this iteration has $\rri=\aleph_1<\ger c$.  To complete the
  proof of the corollary, we shall show that this model has
  $\non N=\ger c$.

  Consider, in the final model, any infinite set $A$ of reals with
  cardinality $<\ger c$; our goal is to prove that $A$ has measure
  zero.  For this purpose, we need to invoke several well-known facts.

First, a finite-support iteration, of length $<\ger c^+$, of
$\sigma$-centered forcing notions is $\sigma$-centered.  This implies
that the standard partial order for forcing \masc\ is
$\sigma$-centered, and therefore that our whole iteration is
$\sigma$-centered.

Second, the regular-open Boolean completion of any $\sigma$-centered
partial order is $\sigma$-centered, and so are all its Boolean
subalgebras.

Third, random real forcing is not $\sigma$-centered.

Combining these three facts, we find that no reals in our final model
are random over the ground model.  This is, in particular, the case
for the reals in $A$.  So, for each $a\in A$, there is a measure-zero
Borel set $N_a$ in the ground model such that the canonical extension
$\tilde N_a$ with the same Borel code contains $a$.

Thanks to the countable chain condition, there is, in the ground
model, a collection \scr C of at most $|A|$ Borel sets, each of
measure zero, such that all of the $N_a$'s are elements of \scr C.
Because the ground model satisfies MA and because $|\scr C|<\ger c$,
the ground model has a measure-zero Borel set $N$ that includes all
the sets from \scr C and, in particular, all the $N_a$'s.

For each $a\in A$, the fact that $N_a\subseteq N$ is preserved when we
pass to the canonical extensions with the same Borel codes in the
final model; that is, $\tilde N_a\subseteq\tilde N$.  In particular,
$A\subseteq\tilde N$.  But $\tilde N$ has, like $N$, measure zero.
This completes the proof that $A$ has measure zero.
\end{pf}

The proof of Lemma~\ref{ccc1} can be modified to show that $\bbb P_I$ is $\sigma$-linked (the details are left as an exercise). But we do not know whether there is a $\sigma$-centered forcing accomplishing the same goal as $\bbb P_I$. For this reason, we do not know whether the analogue of Corollary~\ref{sigmacent:cor} can be proved for $\rrf$ in place of $\rri$, and we leave it as an open question whether it is consistent to have $\rrf < \non N$.

\section{Almost Disjoint Signs} \label{ad}

A natural variant of the rearrangement numbers asks for the minimum
cardinality of a set $C$ of permutations of \bbb N such that, for any
real series $\sum_na_n$, if it has a convergent rearrangement (i.e.,
if it is pcc as defined in Definition~\ref{pcc:df}), then
$\sum_na_{p(n)}$ converges for some $p\in C$.  We show in this section
that this cardinal, unlike the rearrangement  numbers, is provably
equal to the cardinality of the continuum.  That is an immediate
consequence of the following theorem.

\begin{thm}     \label{ad:thm}
  There is a family \scr S of series such that $|\scr S|=\ger c$, that
  each series in \scr S has a convergent rearrangement, but that no
  permutation of \bbb N makes more than one series from \scr S
  converge.
\end{thm}

\begin{pf}
  We use the well-known set-theoretic result that there is a family
  \scr A of \ger c infinite subsets of \bbb N such that any two
  distinct sets from \scr A have finite intersection.  One of several
  easy constructions of such a family proceeds as follows. Instead of
  looking for subsets of \bbb N, we'll get subsets of \bbb Q; they can
  be transferred to \bbb N by any bijection between \bbb N and \bbb Q.
  For each real number $r$, pick a sequence of distinct rationals
  converging to $r$, and let $A_r$ be the range of that sequence.
  Then $\scr A=\{A_r:r\in\bbb R\}$ is as desired.

  Fix a family \scr A as above.  In addition to the fact that distinct
  sets $X,Y\in\scr A$ have $X\cap Y$ finite, we shall need that they
  have $\bbb N\setminus(X\cup Y)$ infinite.  This is easily seen by
  considering a third element $Z$ of $\scr A$ and noting that all but
  finitely many of its elements must be outside $X\cup Y$.

  With these preliminaries out of the way, we proceed to the
  construction of the series required in the theorem.

For each positive integer $i$, let
  $I_i$ be the interval $[2^i+1,2^{i+1}]$, and recall that
  $\sum_{n\in I_i}1/n\geq1/2$.  For each subset $X$ of \bbb N, define a
  series $\sum_na^X_n$ by
\[
a^X_n=
\begin{cases}
  -1/n&\text{if }n\in I_i\text{ for some }i\in X\\
+1/n&\text{otherwise}.
\end{cases}
\]
Consider the series $\sum_na^X_n$ for $X\in\scr A$.  As both $X$ and
its complement $\bbb N\setminus X$ are infinite, the series
$\sum_na^X_n$ includes infinitely many blocks of negative terms of the
form $\sum_{n\in I_i}(-1/n)$, each with sum $\leq-1/2$, and infinitely
many blocks of positive terms of the form $\sum_{n\in I_i}1/n$, each
with sum $\geq1/2$. So the positive and negative parts both diverge,
while the individual terms approach zero,
and therefore the series has a conditionally convergent
rearrangement.

On the other hand, when $X\neq Y$ are distinct elements of \scr A,
then, since $X\cap Y$ is finite, the series $\sum_n(a^X_n+a^Y_n)$ has
only finitely many negative terms and infinitely many blocks of
positive terms (because $\bbb N\setminus(X\cup Y)$ is infinite) with
sum $\geq1/2$ in each block. So it diverges to $+\infty$ under all
permutations.  Therefore, the set $\scr S$ of series $\sum_na^X_n$ for
$X\in\scr A$ is as required in the theorem.
\end{pf}

Note that the series $\sum_na^X_n$ obtained in the proof of the
theorem all diverge by oscillation, since they include blocks of terms
with sums $\leq-1/2$ and blocks with sums $\geq1/2$.  Intuitively,
this argument and the proof of Theorem~\ref{mix-o:thm} suggest that it
is easy to go from conditional convergence to oscillatory divergence
but difficult to go in the other direction.

\section{Shuffles}      \label{shuffles}

The proof of Riemann's Rearrangement Theorem uses only rather special
permutations of \bbb N.  Given a conditionally convergent series, one
uses permutations that keep the relative order of the positive terms
unchanged and also keep the relative order of the negative terms
unchanged.  The only time the relative order of two terms in the
series is changed by the permutation is when one is positive and one
is negative. The following definition formalizes this idea.

\begin{df}      \label{shuffle:df}
Let $A$ and $B$ be two infinite, coinfinite subsets of \bbb N.  The
\emph{shuffle} determined by $A$ and $B$ is the permutation $s_{A,B}$
of \bbb N that maps $A$ onto $B$ preserving order and maps
$\bbb N\setminus A$ onto $\bbb N\setminus B$ preserving order.  That is,
\[
s_{A,B}(n)=
\begin{cases}
  k\th\text{ element of }B&\text{ if }n\text{ is the }k\th\text{
    element of }A\\
  k\th\text{ element of }\bbb N\setminus B&\text{ if }n\text{ is the
  }k\th\text{ element of }\bbb N\setminus A
\end{cases}
\]
\end{df}

We shall be concerned only with the special case where we are
considering a conditionally convergent series $\sum_na_n$ and the
shuffles $s_{A,B}$ under consideration have $B$ equal to
$\{n:a_n>0\}$.  Then $A$ will be the set of locations of positive
terms in the series $\sum_na_{s_{A,B}(n)}$, because
\[
n\in A\iff s_{A,B}(n)\in B\iff a_{s_{A,B}(n)}>0.
\]
In fact, we specialize even further, to alternating series where $B$
is the set of even numbers.  In this situation, we abbreviate
$s_{A,B}$ to $s_A$.  It is the permutation that puts the positive
terms of $\sum_na_n$, in order, into the positions in $A$ and puts the
remaining terms, in order, into $\bbb N\setminus A$.

We could define analogs of all our rearrangement numbers using
shuffles rather than arbitrary permutations, but we actually consider
only the analog of $\rrf$.  That is, we ask how many shuffles $s_A$
are needed to give every conditionally convergent, alternating series
a different, finite sum.  The answer is \ger c, as the following
theorem immediately implies.

\begin{thm}
Consider the conditionally convergent, alternating series
\[
S_\alpha=\sum_n(-1)^n\frac1{(n+1)^\alpha}
\]
for exponents $0<\alpha<1$.  No shuffle $s_A$ makes two of these
series converge to new finite sums. More precisely, if
$0<\alpha<\beta<1$ and if
\[
\sum_n(-1)^{s_A(n)}\frac1{(s_A(n)+1)^\beta}
\]
converges to a finite sum larger (resp.\ smaller) than $S_\beta$, then
\[
\sum_n(-1)^{s_A(n)}\frac1{(s_A(n)+1)^\alpha}
\]
diverges to $+\infty$ (resp.\ to $-\infty$).
\end{thm}

\begin{pf}
We prove the part of the theorem with ``larger'' and $+\infty$; the
proof of the other part is entirely analogous.  Let $\Delta$ be a
positive number such that
\[
\sum_n(-1)^{s_A(n)}\frac1{(s_A(n)+1)^\beta}>
\sum_n(-1)^n\frac1{(n+1)^\beta}  +\Delta.
\]
Let $m$ be a large natural number; just how large $m$ should be will
be determined gradually in the following argument.  Let $x$ be the sum
of the first $m$ positive terms and the first $m$ negative terms in
$S_\beta$, i.e.,
\[
x=\sum_{n=0}^{2m-1}\frac{(-1)^n}{(n+1)^\beta}.
\]
We want to compare $x$ with a certain partial sum $y$ of the series
rearranged by $s_A$, namely the partial sum that ends with the same
negative term $-1/(2m)^\beta$.  This partial sum will have the same
$m$ negative terms as the sum defining $x$ (because $s_A$ is a
shuffle), but it may have more or fewer positive terms; say it has
$m+E(m)$ positive terms, where $E(m)$, the number of excess terms,
might be positive or negative.  These positive terms are, again because
$s_A$ is a shuffle, the first $m+E(m)$ positive terms of the series
$S_\beta$.

For sufficiently large $m$, $x$ will be very close to $S_\beta$
and $y$ will be close to the sum of the series rearranged by $s_A$,
which is more than $S_\beta+\Delta$.  Therefore, taking $m$ large
enough, we have $y>x+\Delta$.  In particular, $E(m)$ must be positive,
i.e., there must be more positive terms in $y$ than in $x$.  Our next
step is to estimate from below the asymptotic size of $E(m)$ for large
$m$.

The excess terms counted by $E(m)$ begin after the $m\th$ positive
term in the original series $S_\beta$, so they are no larger than
$1/(2m)^\beta$.  The sum of these $E(m)$ terms must be more than
$\Delta$, so we have
\[
\frac{E(m)}{(2m)^\beta}>\Delta \quad \text{and so}\quad
E(m)>c\cdot m^\beta
\]
for a suitable positive constant $c$.

With this estimate available, we turn to the other series $S_\alpha$
and its rearrangement by the same $s_A$.  We consider large $m$ and
the partial sums $x$ and $y$ defined as before but with $\alpha$ in
place of $\beta$.  (We can safely use
the same symbols $x$ and $y$ in this new context, as we shall have no
further use for their old meanings.)  It is important to observe that,
because we are using the same shuffle $s_A$ as before, the excess
$E(m)$ is also the same as before, and in particular it obeys the
asymptotic lower bound obtained above, a constant times $m^\beta$.

We use this lower bound to estimate the new $y-x$.  There are two
cases to consider, depending on whether $E(m)\leq m$ or not.
(Actually, the ``not'' case is impossible, but it can be handled
directly just as easily as it can be proved impossible.)

Consider first the case that $E(m)\leq m$. Then the excess terms
counted by $E(m)$ begin with the $m\th$ positive term of $S_\alpha$
and end before the $(2m)\th$ one. In particular, each of these terms
is larger than $1/(4m)^\alpha$. The sum of the excess terms is
therefore asymptotically $\geq c\cdot m^\beta/(4m)^\alpha$, which
tends to $+\infty$ with $m$ because $\alpha<\beta$.  So $y$ grows
without  bound as $m$ increases, which means that the series
$\sum_n(-1)^{s_A(n)}/(s_A(n)+1)^\alpha$ diverges to $+\infty$ as
required.

There remains the case that $E(m)>m$. In this case, instead of adding
all $E(m)$ of the excess terms in $y$, we obtain a lower bound by
adding only the first $m$ of them.  These are, as above, greater than
$1/(4m)^\alpha$, so their sum is at least $m/(4m)^\alpha$, which tends
to $+\infty$ with $m$ because $\alpha<1$.  As in the previous case,
this allows us to conclude the required divergence to $+\infty$.
\end{pf}

\section{Questions}

To conclude the paper, we list some questions that remain open.

\begin{que}
  Is it consistent with ZFC that $\rr_{fi}$, \rrf, and \rri\ are
  different?
\end{que}

Recall that \rr\ can consistently be strictly smaller than $\rr_{fi}$
and \emph{a fortiori} smaller than \rrf\ and \rri, by
Corollary~\ref{separaterr:cor}. But the latter three cardinals are not
separated by any of our results; as far as we know, all three might be
provably equal.  Note that, in this case, the forcing constructions in
Sections~\ref{force1} and \ref{force2} would each achieve the other's
goal (as well as its own).

\begin{que}
  Does ZFC prove that $\rr=\non M$?
\end{que}

\begin{que}
  More generally, are any of the rearrangement  numbers provably equal
  to any previously studied cardinal characteristics of the continuum?
\end{que}

\begin{que}
  Yet more generally, are there provable inequalities between the
  rearrangement numbers and previously studied characteristics, beyond
  those that follow from our results and previously known cardinal
  characteristic inequalities?
\end{que}

\begin{que}
  In particular, are any previously studied characteristics provably
  $\geq\rr_{fi}$?
\end{que}

\end{document}